\magnification=\magstep1
\hsize=16truecm
 
\input amstex
\TagsOnRight
\parindent=20pt
\parskip=2pt plus 1pt
\define\({\left(}
\define\){\right)}
\define\[{\left[}
\define\]{\right]}
\define\e{\varepsilon}
\define\oo{\omega}
\define\const{\text{\rm const.}\,}
\define\supp {\sup\limits}
\define\inff{\inf\limits}
\define\summ{\sum\limits}
\define\prodd{\prod\limits}

\define\bigcupp{\bigcup\limits}

\centerline{\bf
An estimate on the maximum of a nice class of stochastic integrals.}
\smallskip
\centerline{\it P\'eter Major}
\centerline{Alfr\'ed R\'enyi Mathematical Institute of the Hungarian
Academy of Sciences}
\centerline{Budapest, P.O.B. 127 H--1364, Hungary, e-mail:
major\@renyi.hu}
\medskip
{\narrower \noindent {\it Summary:}\/
Let a sequence of iid. random variables $\xi_1,\dots,\xi_n$ be
given on a space $(X,\Cal X)$ with distribution $\mu$ together with
a nice class $\Cal F$ of functions $f(x_1,\dots,x_k)$ of $k$ variables
on the product space  $(X^k,\Cal X^k)$. For all $f\in \Cal F$ we
consider the random integral  $J_{n,k}(f)$ of the function
$f$ with respect to the $k$-fold product of the normalized signed
measure $\sqrt n(\mu_n-\mu)$, where $\mu_n$ denotes the empirical
measure defined by the random variables $\xi_1,\dots,\xi_n$ and
investigate the probabilities $P\(\supp_{f\in \Cal F}|J_{n,k}(f)|>x\)$
for all $x>0$. We show that for nice classes of functions, for
instance if $\Cal F$ is a Vapnik--\v{C}ervonenkis class, an almost
as good bound can be given for these probabilities as in the case when
only the random integral of one function is considered. \par}
 
\beginsection 1. Introduction. Formulation of the main results
 
The following problem is studied in this paper: Let a probability
measure $\mu$ be given on a measure space $(X,\Cal X)$, take a
sequence $\xi_1,\dots,\xi_n$ of independent, identically distributed
$(X,\Cal X)$ valued random variables with distribution~$\mu$, and
define the empirical measure $\mu_n$,
$$
\mu_n(A)=\dfrac1n\#\{j\:\,\xi_j\in A,\;1\le j\le n\},\quad A\in\Cal X,
$$
of the sample $\xi_1,\dots,\xi_n$. Let us take a nice set $\Cal F$
of measurable functions $f(x_1,\dots,x_k)$ on the $k$-fold product
space $(X^k,\Cal X^k)$ and define the integrals $J_{n,k}(f)$ of
the functions $f\in \Cal F$ with respect to the $k$-fold product
of the normalized empirical measure $\mu_n$ by the formula
$$
\align
J_{n,k}(f)&=\dfrac{n^{k/2}}{k!} \int'
f(u_1,\dots,u_k)(\mu_n(\,du_1)-\mu(\,du_1))\dots
(\mu_n(\,du_k)-\mu(\,du_k)),\\
&\qquad\text{where the prime in $\tsize\int'$ means that the
diagonals } u_j=u_l,\; 1\le j<l\le k,\\
&\qquad\text{are omitted from the domain of integration.} \tag1.1
\endalign
$$
In this work I try to give a good estimate on the probabilities
$P\(\supp_{f\in\Cal F}|J_{n,k}(f)|>x\)$ for all $x>0$. To formulate
the main result of the paper first I introduce the following definition.
\medskip\noindent
{\bf Definition of $L_p$-dense classes of functions.} {\it Let us
have a measure space $(Y,\Cal Y)$ and a set $\Cal G$ of $\Cal Y$
measurable functions on this space. We call $\Cal G$ an $L_p$-dense
class with parameter $D$ and exponent $L$ if for all numbers
$1\ge\e>0$ and probability measures $\nu$ on the space $(Y,\Cal Y)$
there exists a finite $\e$-dense subset $\Cal G_{\e,\nu}=\{g_1,\dots,
g_m\} \subset \Cal G$ in the space $L_p(Y,\Cal Y,\nu)$ consisting of
$m\le D\e^{-L}$ elements, i.e. such a set $\Cal G_{\e,\nu}\subset \Cal
G$ for which $\inff_{g_j\in \Cal G_{\e,\nu}}\int |g-g_j|^p\,d\nu<\e^p$
for all functions $g\in \Cal G$. (Here the set $\Cal G_{\e,\nu}$ may
depend on the measure $\nu$, but its cardinality is bounded by a
number depending only on $\e$.)}
\medskip
In this paper we shall work with such classes of functions $\Cal F$
which contain only functions with absolute value less than or equal
to~1. In this case $\Cal F$ is an $L_p$-dense class of functions for
all $1\le p<\infty$ (with an exponent and a parameter depending
on~$p$) if there is a number $1\le p<\infty$ for which it is
$L_p$-dense. We shall formulate our statements mainly for $L_p$-dense
classes of functions with the parameter $p=2$, since this seems to be
the most convenient choice. Our main result is the following
\medskip\noindent
{\bf Theorem.} {\it Let us have a non-atomic measure $\mu$ on the
space $(X,\Cal X)$ together with an $L_2$-dense class $\Cal F$ of
functions $f=f(x_1,\dots,x_k)$ of $k$ variables with some parameter
$D$ and exponent $L$ on the product space $(X^k,\Cal X^k)$ which
consists of at most countably many functions, and satisfies the
conditions
$$
\|f\|_\infty=\supp_{x_j\in X,\;1\le j\le k}|f(x_1,\dots,x_k)|\le 1,
\qquad \text{for all } f\in \Cal F \tag1.2
$$
and
$$
\|f\|_2^2=Ef^2(\xi_1,\dots,\xi_k)=\int f^2(x_1,\dots,x_k)
\mu(\,dx_1)\dots\mu(\,dx_k)\le \sigma^2 \qquad \text{for all }
f\in \Cal F \tag1.3
$$
with some constant $\sigma>0$. Let us also assume that the parameter
$D$ of the $L_2$-dense class $\Cal F$ satisfies the condition
$$
D\le n^\beta \quad\text{with some } \beta\ge0.  \tag1.4
$$
Then there exist some constants $C=C(k)>0$, $\alpha=\alpha(k)>0$
and $M=M(k)>0$ depending only on the parameter $k$ such that the
supremum of the random integrals $J_{n,k}(f)$, $f\in \Cal F$,
defined by formula (1.1) satisfies the inequality
$$
\aligned
P\(\supp_{f\in\Cal F}|J_{n,k}(f)|\ge x\)&\le CD \exp\left\{-\alpha
\(\frac x{\sigma}\)^{2/k}\right\} \\
&\qquad \text{if}\quad n\sigma^2\ge
\(\frac x\sigma\)^{2/k} \ge M(L+\beta+1)^{3/2}\log\frac2\sigma,
\endaligned \tag1.5
$$
where $\beta$ is the number in (1.4), and the number $D$ in
formula~(1.5) agrees with the parameter of the $L_2$-dense
class~$\Cal F$.}
\medskip
The condition that $\Cal F$ is a countable class of functions can be
weakened. To formulate such a result the following definition will be
introduced.
\medskip\noindent
{\bf Definition of countable majorizability.} {\it A class of
functions $\Cal F$ is countably majorizable in the space
$(X^k,\Cal X^k,\mu^k)$ if there exists a countable subset $\Cal
F'\subset \Cal F$ such that for all numbers $x>0$ the sets
$A(x)=\{\oo\:\supp_{f\in \Cal F}|J_{n,k}(f)(\oo)|\ge x\}$ and
$B(x)=\{\oo\:\supp_{f\in \Cal F'}|J_{n,k}(f)(\oo)|\ge x\}$ satisfy
the identity $P(A(x)\setminus B(x))=0$.}
\medskip
Clearly, $B(x)\subset A(x)$. In the above definition we demanded
that for all $x>0$ the set $B(x)$ is almost as large as $A(x)$.
Now the following corollary of the Theorem is given.
\medskip\noindent
{\bf Corollary 1 of the Theorem.} {\it Let a class of functions
$\Cal F$ satisfy the conditions of the Theorem with the only
exception that instead of the condition about the countable
cardinality of $\Cal F$ it is assumed that $\Cal F$ is countably
majorizable in the space $(X^k,\Cal X^k,\mu^k)$. Then $\Cal F$
satisfies the Theorem.}
\medskip
 
The condition that the class of functions $\Cal F$ is countable was
imposed to avoid some unpleasant measure theoretical difficulties
which would arise if we had to work with possibly non-measurable
sets. On the other hand, I have the impression that Corollary~1 can
be applied in all investigations where an estimate about the
supremum of multiple random integrals with respect to a normalized
empirical measure is needed. It is not difficult to prove that
Corollary~1 follows from the Theorem. To do this we have to show
that if $\Cal F$ is an $L_2$-dense class with some parameter $D$
and exponent $L$, and $\Cal F'\subset \Cal F$, then $\Cal F'$ is
also an $L_2$-dense class with the same exponent $L$, only with a
possibly different parameter~$D'$.
 
To prove this statement let us choose for all numbers $1\ge\e>0$
and probability measures $\nu$ on $(Y,\Cal Y)$ some functions
$f_1,\dots,f_m\in \Cal F$ with $m\le D\(\frac\e2\)^{-L}$ elements,
such that the sets $\Cal D_j=\left\{f\:\int |f-f_j|^2\,d\nu\le
\(\frac\e2\)^2\right\}$ satisfy the relation $\bigcupp_{j=1}^m \Cal
D_j=Y$. For all sets $\Cal D_j$ for which $\Cal D_j\cap \Cal F'$ is
non-empty choose a function $f'_j\in \Cal D_j\cap \Cal F'$. In such a
way we get a collection of functions $f'_j$ from the class $\Cal F'$
containing at most $2^LD\e^{-L}$ elements which satisfies the
condition imposed for $L_2$-dense classes with exponent $L$ and
parameter $2^LD$ for this number $\e$ and measure $\nu$.
\medskip
 
The following Corollary of the Theorem may be of special interest.
It is similar to some results of paper [2] or Theorem~5.3.14 in~[4].
\medskip\noindent
{\bf Corollary 2 of the Theorem.} {\it Let us consider a non-atomic
probability measure $\mu$ on a measure space $(X,\Cal X)$ and an
$L_2$-dense class $\Cal F$ of functions on the $k$-fold product
space $(X^k,\Cal X^k)$ with some exponent $L$ and parameter $D$
which is either countable or countably majorizable. Let us also
assume that $\supp_{x_j\in X,\;1\le j\le k}|f(x_1,\dots,x_k)|\le 1$
for all $f\in \Cal F$. Then the supremum of the random stochastic
integrals $J_{n,k}(f)$, $f\in\Cal F$, satisfies the inequality
$$
P\(\supp_{f\in\Cal F}|J_{n,k}(f)|\ge x\)\le C e^{-\alpha x^{2/k}}
\tag1.6
$$
for all $x>0$ with some constants $\alpha=\alpha(k)>0$ and
$C=C(k,D,L)$ depending on the parameter $k$ on the exponent and
parameter of the $L_2$-dense class $\Cal F$.}
\medskip\noindent
{\it Proof of Corollary 2.} Let us first assume that $D\le n$, and
apply the result of the Theorem with $\sigma=1$. Then conditions (1.3)
and (1,4) of the Theorem hold. Also the first part of the condition
in (1.5) holds if $x\le n^{k/2}$, and $P(|J_{n,k}(f)|>x)=0$ if
$x>n^{k/2}$. The second condition of (1.5) is satisfied if $x^{2/k}
\ge M((L+\beta+1)^{3/2}\log 2$, hence relation (1.6) holds if
$x\ge\const$ with an appropriate constant. If the number $C$ in (1.6)
is chosen sufficiently large, then the right-hand side of (1.6) is
greater than 1 for $x\le \const$ In the case $n\le D$ the random
integral $|J_{n,k}(f)|$ is less than $\frac{2^kD^{k/2}}{k!}$. Hence
the statement of Corollary~2 holds for all $x>0$ with an appropriate
choice of the parameter~$C$. \medskip
 
In the Theorem we have considered the supremum of multiple
random integrals for a nice class of functions of $k$ variables
with respect to the $k$-fold product of a normalized empirical
measure. It was shown that if the variances of the random integrals
we have considered are less than some number $\sigma^2>0$, then under
some additional conditions this supremum takes a value larger than
$x$ with a probability less than $P(C\sigma\eta>x)$, where $\eta$ is
a standard normal random variable, and $C=C(k)$ is a universal
constant depending only on the multiplicity $k$ of the random
integrals. This is the sharpest estimate we can expect. Moreover,
this estimate seems to be sharp also in that respect that the
conditions imposed for its validity cannot be considerably
weakened. If condition (1.2) does not hold or $n\sigma^2<
\(\frac x\sigma\)^{2/k}$, then the estimate of the Theorem may
not hold any longer even if the class of functions $\Cal F$ contains
only one function. In such cases there exist examples for which the
probability $P(J_{n,k}(f)>x)$ is too large. Indeed, in such cases
it may happen that the value of relatively few members of the
sample take the random integral larger than~$x$ with relatively
large probability, and the remaining part of the sample does not
diminish it. Here I do not work out the details of such examples.
 
If $\(\frac x\sigma\)^{2/k}<M\log\frac2\sigma$ with a not too large
number $M>0$, then the estimate of the Theorem may be violated again,
but in this case the reason for it is that the supremum of many small
random variables may be large. To understand this let us consider the
following analogous problem. Take a Wiener process $W(t)$, $0\le t\le1$,
and consider the supremum of the expressions $W(t)-W(s)=\int
f_{s,t}(u)W(\,du)=\bar J(f_{s,t})$, with the functions $f_{s,t}(\cdot)$
on the interval $[0,1]$ defined by the formula $f_{s,t}(u)=1$ if $s\le
u\le t$, $f_{s,t}(u)=0$ if $0\le u<s$ or $t<u\le1$. If we consider the
class of functions $\Cal F_\sigma=\{f_{s,t}\:\int
f^2_{s,t}(u)\,du=t-s\le\sigma^2\}$, then it is natural to expect that
$P\(\supp_{f_{s,t}\in\Cal F_\sigma} \bar J(f_{s,t})>x\)\le e^{-\const
(x/\sigma)^2}$. However, this relation does not hold if
$x=x(\sigma)<(1-\e)\sqrt{2\log\frac1\sigma}\sigma$ with some $\e>0$.
In such cases $P\(\supp_{f_{s,t}\in\Cal F_\sigma}\bar
J(f_{s,t})>x\)\to1$, as $\sigma\to0$. This can be proved relatively
simply with the help of the estimate $P(\bar J(f_{s,t})>x(\sigma))
\ge\const \sigma^{1-\e}$ if $|t-s|=\sigma^2$ and the independence of
the random integrals $\bar J(f_{s,t})$ if the functions $f_{s,t}$ are
indexed by such pairs $(s,t)$ for which the intervals $(s,t)$ are
disjoint.
 
Some additional work would show that a similar picture arises if we
integrate with respect to the normalized empirical measure
of a sample with uniform distribution on the interval $[0,1]$ instead
of a Wiener process. This yields an example for an $L_2$-dense
class of functions in the case $k=1$ for which the estimate of the
Theorem does not hold any longer if $\(\frac
x\sigma\)^{2/k}<M\log\frac2\sigma$ with some $M<\sqrt2$.
At a heuristic level it is clear that such an example can be given
also for $k>1$, and the number $M$ in condition (1.5) has to be
chosen larger if we want that the Theorem hold also for an
$L_2$-dense class of functions $\Cal F$ with a large exponent~$L$.
(In this paper I did not try to find the best possible condition of the
 Theorem in the right-hand side inequality of~(1.5).)
 
One would like to see some interesting examples when the Theorem is
applicable and to have some methods to check the conditions of the
Theorem. It is useful to know that if $\Cal F$ is a
Vapnik--\v{C}ervonenkis class of functions whose absolute values
are bounded by 1, then $\Cal F$ is an $L_2$-dense class.
 
To formulate the above statement more explicitly let us recall that
a class of subsets $\Cal D$ of a set $S$ is a Vapnik--\v{C}ervonenkis
class if there exist some constants $B>0$ and $K>0$ such that for all
integers $n$ and sets $S_0(n)=\{x_1,\dots,x_n\}\subset S$ of
cardinality $n$ the collection of sets of the form $S_0(n)\cap D$,
$D\in\Cal D$, contains no more than $Bn^K$ subsets of $S_0(n)$. A
class of real valued functions $\Cal F$ on a space $(Y,\Cal Y)$
is a Vapnik--\v{C}ervonenkis class if the graphs of these functions
is a Vapnik--\v{C}ervonenkis class, i.e.\ if the sets $A(f)=\{(y,t)\:
y\in Y,\;\min(0,f(y))\le t\le\max(0,f(y))\}$, $f\in \Cal F$, constitute
a Vapnik--\v{C}ervonenkis class of sets on the product space $Y\times
R^1$.
 
An important result of Dudley states that a Vapnik--\v{C}ervonenkis
class of functions whose absolute values are bounded by 1 is an
$L_1$-dense class. The parameter and exponent of this $L_1$-dense
class can be bounded by means of the constants $B$ and $K$ appearing
in the definition of Vapnik--\v{C}ervonenkis classes. On the other
hand, an $L_1$-dense class of functions bounded by 1 is also an
$L_2$-dense class (with possibly different exponent and parameter),
since $\int|f-g|^2\,d\nu\le2\int|f-g|\,d\nu$ in this case. Dudley's
result, whose proof can be found e.g.\ in Chapter~II of Pollard's
book~[9] (the 25$^\circ$ approximation lemma contains this result in
a slightly more general form) is useful for us, because there are
results which enable us to prove that certain classes of functions
constitute a Vapnik--\v{C}ervonenkis class.
 
\medskip
 
This work is a continuation of my paper~[8], where this question
was discussed in the special case when $\Cal F$ contains only one
function. Here Theorem~1 (or its equivalent version Theorem~$1'$)
of [8] will be applied, but no additional argument of that work is
needed. As I have mentioned in~[8], the investigation of this paper
was motivated by some non-parametric maximum likelihood estimate
problems. Earlier I could only prove a much weaker version of this
result in [7].
 
I found some results similar to that of this paper
in the work of Arcones and Gine~[2], where the tail-behaviour of
the supremum of degenerated $U$-statistics was investigated if the
kernel functions of these $U$-statistics constitute a
Vapnik--\v{C}ervonenkis class. But the bounds of that paper do not
give a better estimate if we have the additional information that the
variances of the $U$-statistics we consider are small. On the other
hand, one of the main goal of the present paper was to prove such
estimates. (Let me remark that formula~(1.3) imposes a condition on
the variances of the random integrals we consider in this paper. See
Lemma~3 in~[8].) I know of one work where the dependence of the
estimate on the variance was investigated in a similar case. This is
Alexander's paper~[1], where the problem of the present paper was
studied in the special case $k=1$. Alexander proved in this
case a sharper result. He also studied the case of non-identically
distributed random variables and gave an upper bound for the
distribution function of the supremum of random integrals with
almost as good constants as in the case of a single random integral.
Probably a similar result also holds for multiple stochastic
integrals, but the proof requires a more careful analysis. Alexander's
paper was interesting for me first of all, because I learned some ideas
from it which I strongly needed in the present work. On the other hand,
I also needed some new arguments, because in the study of multiple
stochastic integrals some new difficulties had to be overcome.
 
This paper consists of six sections and an appendix. In Section~2
the Theorem is reduced to a simpler statement formulated in
Proposition~3. Section~3 contains some important results needed
in the proof, and the main ideas of the proof are explained
there. In particular, the proof of Proposition~3 is reduced to
another statement formulated in Proposition~4. Proposition~4
is proved simultaneously with another result described in
Proposition~5. To make the proof more transparent first I give
it in the special case $k=1$ in Section~4. Sections~5 and~6 contain
the proof of Propositions~4 and~5 in the general case. In Section~5
it is shown how a symmetrization argument can be applied to prove
these results, and finally the proof is completed in Section~6. The
Appendix contains the proof of an estimate about the tail behaviour of
the distribution of homogeneous polynomials of Rademacher functions.
 
\beginsection
2. Reduction of the Theorem to a simpler result
 
I shall prove with the help of a natural argument, called the Chaining
argument in the literature, and the result Theorem~$1'$ in paper~[9]
the following result.
\medskip\noindent
{\bf Proposition 1.} {\it Let us fix some number $\bar A\ge2^k$,
and assume that a class of functions $\Cal F$ satisfies the
conditions of the Theorem with a number $M$ in these conditions
which may depend also on $\bar A$. Then a number $0\le\bar\sigma\le
\sigma\le1$ and a collection of functions $\Cal F_{\bar\sigma}
=\{f_1,\dots,f_m\}\subset \Cal F$ with $m\le D\bar\sigma^{-L}$
elements can be chosen in such a way that the sets
$\Cal D_j=\{f\:f\in \Cal F,\int|f-f_j|^2\,d\mu\le\bar\sigma^2\}$,
$1\le j\le m$, satisfy the relation $\bigcupp_{j=1}^m\Cal D_j=\Cal F$,
and
$$
\aligned
P&\(\sup_{f\in\Cal F_{\bar\sigma}} |J_{n,k}(f)|\ge \frac x{\bar
A}\)\le 2CD\exp\left\{-\alpha\(\frac x{4\bar
A\sigma}\)^{2/k}\right\}  \\
&\qquad \qquad \text{if}\quad n\sigma^2\ge \(\frac x\sigma\)^{2/k}
\ge ML\log\frac2\sigma
\endaligned \tag2.1
$$
with the constants $\alpha=\alpha(k)$, $C=C(k)$ appearing in
Theorem~$1'$ of~[8]  and the exponent $L$ and parameter $D$ of the
$L_2$-dense class $\Cal F$ if the constant $M=M(k,\bar A)$ is
chosen sufficiently large. Beside this, also the inequalities
$64\(\frac x{\bar A\bar\sigma}\)^{2/k}\ge n\bar\sigma^2\ge
\(\frac x{\bar A\sigma}\)^{2/k}$ and
$n\bar\sigma^2\ge\frac{M^{2/3}(L+\beta+1)\log n}{50\bar A^{4/3}}$
hold, provided that $n\sigma^2\ge \(\frac x\sigma\)^{2/k}\ge
M(L+\beta+1)^{3/2}\log\frac2\sigma$.}
\medskip\noindent
{\it Remark:}\/ The introduction of the number $\bar A\ge2$ in
Proposition~1 may seem a bit artificial. Its role is to guarantee
that such a number $\bar\sigma$ could be defined in Proposition~1
which satisfies the inequality $\(\frac x{\bar\sigma}\)^{2/k}\ge A
n\bar\sigma^2$ with a sufficiently large previously fixed constant
$A=A(k)$.
\medskip\noindent
{\it Proof of Proposition 1.} For all $p=0,1,2,\dots$ choose a set
$\Cal F_p=\{f_{p,1},\dots,f_{p,m_p}\}\subset\Cal F$ with $m_p\le D\,
4^{pL}\sigma^{-L}$ elements in such a way that $\inff_{1\le j\le m_p}
\int (f-f_{p,j})^2\,d\mu\le 16^{-p}\sigma^2$ for all $f\in\Cal F$.
For all pairs $(j,p)$, $p=1,2,\dots$, $1\le j\le m_p$, choose a
precedor $(j',p-1)$, $j'=j'(j,p)$, $1\le j'\le m_{p-1}$, in such a
way that the functions $f_{j,p}$ and $f_{j',p-1}$ satisfy the
relation $\int|f_{j,p}-f_{j',p-1}|^2\,d\mu\le \sigma^2 16^{-p}$.
Then we have $\int\(\frac{f_{j,p}-f_{j',p-1}}2\)^2\,d\mu\le\frac14
\sigma^2 16^{-p}$ and $\supp_{x_j\in X,\,1\le j\le k}\left|
\frac{f_{j,p}(x_1,\dots,x_k)-f_{j',p-1}(x_1,\dots,x_k)}2\right|\le 1$.
Theorem~$1'$ of [8] yields that
$$
\align
P(A(j,p))&=P\(|J_{n,k}(f_{j,p}-f_{j',p-1})|\ge \frac{2^{-(1+p)}x}
{\bar A}\) \le C \exp\left\{-\alpha\(\frac{2^{p-1}x}{\bar A
\sigma}\)^{2/k} \right\}\\
&\qquad \text {if}\quad \frac{n\sigma^2 2^{-4p}}4\ge \(\frac {2^{p-1}x}
{\bar A\sigma}\)^{2/k},\quad 1\le j\le m_p,\; p=1,2,\dots,    \tag2.2
\endalign
$$
and
$$
\aligned
P(B(s))&=P\(|J_{n,k}(f_{0,s})|\ge \frac x{2\bar A}\)\le
C\exp\left\{-\alpha\(\frac x{2\bar A\sigma}\)^{2/k}\right\},
\quad 1\le s\le m,\\
&\qquad\qquad\qquad\text{if} \quad n\sigma^2\ge \(\frac x{2\bar
A\sigma}\)^{2/k}. \endaligned\tag2.3
$$
Choose the integer number $R$, $R\ge0$, in such a way that
$2^{(4+{2/k})(R+1)}\(\frac{x}{\bar A\sigma}\)^{2/k} \ge
\frac{n\sigma^2}{2^{2-2/k}}\ge 2^{(4+2/k)R}\(\frac{x}{\bar
A\sigma}\)^{2/k}$, define $\bar\sigma^2=16^{-R}\sigma^2$ and $\Cal
F_{\bar\sigma}=\Cal F_R$. (As $n\sigma^2\ge\(\frac x\sigma\)^{2/k}$
and $\bar A\ge2^k$ by our conditions, there exists such a
non-negative number $R$.) Then the cardinality~$m$ of the set $\Cal
F_{\bar\sigma}$ is clearly not greater than $D\bar\sigma^{-L}$,
and $\bigcupp_{j=1}^m \Cal D_j=\Cal F$. Beside this, the number
$R$ was chosen in such a way that the inequalities
(2.2) and (2.3) can be applied for $1\le p\le R$. Hence the
definition of the precedor of a pair $(j,p)$ implies that
$$
\align
&P\(\sup_{f\in\Cal F_{\bar\sigma}} |J_{n,k}(f)|\ge \frac x{\bar A}\)
\le P\(\bigcup_{p=1}^R\bigcup_{j=1}^{m_p}A(j,p)
\cup\bigcup_{s=1}^mB(s)\) \\
&\qquad \le \sum_{p=1}^R\sum_{j=1}^{m_p}P(A(j,p))+\sum_{s=1}^mP(B(s))
\le \sum_{p=1}^{\infty} CD\,4^{pL}\sigma^{-L}
\exp\left\{-\alpha\(\frac{2^{p-1}x}{\bar A\sigma}\)^{2/k}
\right\}\\
&\qquad\qquad +CD\sigma^{-L}\exp\left\{-\alpha\(\frac
x{2\bar A\sigma}\)^{2/k}\right\}.
\endalign
$$
If the condition $\(\frac x\sigma\)^{2/k}\ge M(L+1)^{3/2}
\log\frac2\sigma$ holds with a sufficiently large constant
$M$ (depending on $\bar A$), then the inequalities
$$
4^{pL}\sigma^{-L}\exp\left\{-\alpha\(\frac{2^{p-1}x}{\bar
A\sigma}\)^{2/k} \right\}\le 2^{-p} \exp\left\{-\alpha\(\frac{2^{p}x}
{4\bar A \sigma}\)^{2/k} \right\}
$$
hold for all $p=1,2,\dots$, and
$$
\sigma^{-L}\exp\left\{-\alpha\(\frac x{2\bar A\sigma}\)^{2/k}\right\}
\le\exp\left\{-\alpha\(\frac x{4\bar A\sigma}\)^{2/k}\right\}.
$$
Hence the previous estimate implies that
$$
\align
&P\(\sup_{f\in\Cal F_{\bar\sigma}} |J_{n,k}(f)|\ge \frac x{\bar A}\)
\le\sum_{p=1}^{\infty}CD 2^{-p}
\exp\left\{-\alpha\(\frac{2^{p}x}{4\bar A \sigma}\)^{2/k}
\right\}\\
&\qquad +CD\exp\left\{-\alpha\(\frac x{4\bar A
\sigma}\)^{2/k}\right\} \le 2CD \exp\left\{-\alpha\(\frac x{4
\bar A\sigma}\)^{2/k}\right\},
\endalign
$$
and relation (2.1) holds. We have
$$
\align
n\bar\sigma^2&=2^{-4R} n\sigma^2\le
2^{-4R}\cdot2^{(4+2/k)(R+1)+2-2/k}\(\frac{x}{\bar A\sigma}\)^{2/k}=
2^6\cdot 2^{2R/k}\(\frac{x}{\bar A \sigma}\)^{2/k}\\
&=2^6\cdot \(\frac\sigma{\bar\sigma}\)^{1/k}\(\frac{x}{\bar A
\sigma}\)^{2/k}=2^6\cdot \(\frac{\bar\sigma}\sigma\)^{1/k}
\(\frac{x}{\bar A \bar\sigma}\)^{2/k},
\endalign
$$
hence $n\bar\sigma^2\le  2^6 \(\frac{x}{\bar A\bar\sigma}\)^{2/k}$.
Beside this,
$$
n\bar\sigma^2=2^{-4R}n\sigma^2\ge
2^{2-2/k}\cdot 2^{-4R}\cdot 2^{4R+2R/k}\(\frac x{\bar
A\sigma}\)^{2/k}\ge \(\frac x{\bar A\sigma}\)^{2/k}.
$$
It remained to show that $n\bar\sigma^2\ge
\frac{M^{2/3}((L+\beta)\log n+1)}{50A^{4/3}}$.
 
This inequality clearly holds under the conditions of Proposition~1
if $\sigma\le n^{-1/3}$, since in this case $\log\frac2\sigma\ge
\frac{\log n}3$, and $n\bar\sigma^2\ge\(\frac x {\bar A\sigma}\)^{2/k}
\ge\bar A^{-2/k} M(L+\beta+1)^{3/2}\log \frac2\sigma\ge
\frac13\bar A^{-2/k} M(L+\beta+1)\log n$. If $\sigma\ge
n^{-1/3}$, then the inequality $2^{(4+2/k)R}\(\frac x{\bar
A\sigma}\)^{2/k} \le \frac{n\sigma^2}{2^{2-2/k}}$ holds. Hence
$2^{-4R}\ge 2^{(2-2/k))/(4+2/k)}  \[\dfrac{\(\frac
x{\bar A\sigma}\)^{2/k}}{n\sigma^2}\]^{4/(4+2/k)}$, and
$$
n\bar\sigma^2=2^{-4R}n\sigma^2\ge\frac1{\bar A^{4/3}}
(n\sigma^2)^{1-\gamma}\[\(\frac x\sigma\)^{2/k}\]^\gamma \text{ with }
\gamma=\frac4{4+\frac2k}\ge\frac23.
$$
Since $n\sigma^2\ge(\frac x\sigma)^{2/k}\ge\frac M3(L+\beta+1)^{3/2}$,
and $n\sigma^2\ge n^{1/3}$, the above estimates yield that
$n\bar\sigma^2\ge \bar A^{-4/3} (n\sigma^2)^{1/3}\[\(\frac
x\sigma\)^{2/k}\]^{2/3}\ge \bar A^{-4/3}n^{1/9}\(\frac M3\)^{2/3}
(L+\beta+1) \ge\frac{M^{2/3}(L+\beta+1)\log n}{50 \bar A^{4/3}}$.
\medskip
Now I formulate Proposition~2 and show that the Theorem follows from
Propositions~1 and~2.
\medskip\noindent
{\bf Proposition 2.}
{\it Let us have a non-atomic measure $\mu$ on the space $(X,\Cal X)$
together with a sequence of independent and $\mu$ distributed random
variables $\xi_1,\dots,\xi_n$ and an $L_2$-dense class of functions
$f=f(x_1,\dots,x_k)$ of $k$ variables with some parameter $D$ and
exponent $L$ on the product space $(X^k,\Cal X^k)$ which consists of
at most countably many functions and satisfies conditions (1.2),
(1.3) and (1.4) with some $\sigma>0$, and
$n\sigma^2>K((L+\beta)\log n+1)$ with a sufficiently large number
$K=K(k)$. Then there exists some number $\gamma=\gamma(k)>0$ and
threshold index $A_0=A_0(k)>0$ depending only on the order~$k$ of
the stochastic integrals we consider such that
$$
P\(\sup_{f\in\Cal F}|J_{n,k}(f)|\ge A n^{k/2}\sigma^{k+1}\)\le
e^{-\gamma A^{1/2k}n\sigma^2}\quad \text{if } A\ge A_0.
$$
} \medskip
 
In the proof of the Theorem we exploit our freedom in the choice
of the parameters in Propositions~1 and~2. Let us choose a number
$\bar A_0$ such that $\bar A_0\ge A_0$ and $\gamma\bar
A_0^{1/2k}\ge\frac1K$ with the numbers $A_0$, $K$ and $\gamma$ in
Proposition~2. Proposition~1 will be applied with such a number
$\bar A$ for which the inequalities $\(\frac x{\bar\sigma}\)^{2/k}
\ge\frac {\bar A^{2/k}}{64}n \bar\sigma^2\ge (4\bar A_0)^{2/k}n\bar
\sigma^2$ hold with the above fixed parameter $\bar A_0$ and the
number $\bar\sigma$ defined in the proof of Proposition~1. (Here
and in the sequel we shall assume that the number $x$ satisfies the
condition $n\sigma^2\ge \(\frac x\sigma\)^{2/k}\ge
M(L+\beta+1)^{3/2}\log\frac2\sigma$ imposed both in Proposition~1
and in the Theorem.) Choose such a number $M$ in Proposition~1
(and as a consequence in the Theorem too) for which also the
inequality $n\bar\sigma^2\ge\frac{M^{2/3}(L+\beta+1)\log n}{50\bar
A^{4/3}}\ge K((L+\beta)\log n+1)$ holds with the number $K$
appearing in the conditions of Proposition~2. Proposition~1 will be
applied with the class of functions $\Cal F$, the numbers $\sigma$
and $M$ considered in the Theorem and a number~$\bar A$ satisfying
the above property while Proposition~2 with the above chosen
number $\bar A_0$, the number $\bar\sigma$ and the sets of
functions $\Cal D_j$ defined in Proposition~1. More precisely, we
apply Proposition~2 for the sets of functions $\frac{g-f_j}2$ where
$g\in \Cal D_j$ and $f_j$ is the `center' of the set $\Cal D_j$
appearing in the definition of the set $\Cal D_j$ in Proposition~1.
Observe that these functions constitute an $L_2$-dense class of
functions with exponent $L$ and parameter~$D$.
 
Since $\(1-\frac 1 {\bar A}\)x\ge\frac x2\ge 2\bar
A_0n^{k/2}\bar\sigma^{k+1}$ Propositions~1 and~2 with the above
parameters yield that
$$
\aligned
P\(\supp_{f\in\Cal F}|J_{n,k}(f)|\ge x\)&\le
P\(\sup_{f\in\Cal F_{\bar\sigma}} |J_{n,k}(f)|\ge \frac x{\bar A}\)\\
&\qquad +\sum_{j=1}^m P\(\sup_{g\in\Cal
D_j}\left|J_{n,k}\(\frac{f_j-g}2\)\right| \ge \bar A_0
n^{k/2}\bar\sigma^{k+1}\) \\
&\le 2CD\exp\left\{-\alpha\(\frac x{4\bar A\sigma}\)^{2/k}\right\}
+D\bar\sigma^{-L} e^{-\gamma\bar A_0^{1/2k}n\bar\sigma^2}.
\endaligned \tag2.4
$$
Let us understand how the second term at the right-hand side of
(2.4) can be estimated. The condition $n\bar\sigma^2\ge
K((L+\beta)\log n+1)$ implies that $\bar\sigma\ge n^{-1/2}$, and
by our choice of $\bar A_0$ we have $\gamma \bar
A_0^{1/2k}n\bar\sigma^2\ge \frac1Kn\bar\sigma^2 \ge L\log n\ge
2L\log\frac1{\bar \sigma}$, i.e. $\bar\sigma^{-L}\le e^{\gamma\bar
A_0^{1/2k}n\bar\sigma^2/2}$. As we have seen in Proposition~1
$n\bar\sigma^2 \ge\(\frac x{\bar A\sigma}\)^{2/k}$. The above
relations imply that $\bar\sigma^{-L} e^{-\gamma\bar A_0^{1/2k}n
\bar\sigma^2}\le e^{-\gamma\bar A_0^{1/2k}n\bar\sigma^2/2}\le
\exp\left\{-\frac\gamma2 \bar A_0^{1/2k} \bar A^{-2/k}\(\frac
x\sigma\)^{2/k}\right\}$. Then relation (2.4) gives that
$$
\align
P\(\supp_{f\in\Cal F}|J_{n,k}(f)|\ge x\)&\le 2CD\exp
\left\{-\frac\alpha{(4\bar A)^{2/k}}\(\frac x\sigma\)^{2/k}\right\}
\\ &\qquad+D\exp\left\{-\frac\gamma2\bar A_0^{1/2k}\bar A^{-2/k}
\(\frac x\sigma\)^{2/k}\right\}.
\endalign
$$
 
The last formula means that under the conditions of the Theorem
formula (1.5) holds (with some new appropriately defined constant
$\alpha>0$), and this is what we had to prove.
\medskip
 
It remained to prove Proposition~2. Its proof requires some new
ideas, and the remaining part of the paper deals with this problem.
There is a  counterpart of this result about so-called degenerate
$U$-statistics. The study of degenerate $U$-statistics is technically
simpler. Hence I formulate this result about $U$-statistics in
Proposition~3 and show that it implies Proposition~2.
 
First I recall some notions we need to formulate Proposition~3.
Let us have a sequence of independent and identically distributed
random variables $\xi_1,\xi_2,\dots$ with distribution~$\mu$ on a
measurable space $(X,\Cal X)$ together with a function
$f=f(x_1,\dots,x_k)$ on the $k$-th power $(X^k,\Cal X^k)$ of the
space $(X,\Cal X)$. We define with their help the
$U$-statistic~$I_{n,k}(f)$ of order~$k$, as
$$
I_{n,k}(f)=\frac1{k!}\summ\Sb 1\le j_s\le n,\; s=1,\dots, k\\
j_s\neq j_{s'} \text{ if } s\neq s'\endSb
f\(\xi_{j_1},\dots,\xi_{j_k}\).   \tag2.5
$$
(The function $f$ in this formula will be called the kernel
function of the $U$-statistic.)
 
A real valued function $f=f(x_1,\dots,x_k)$ on the $k$-th power
$(X^k,\Cal X^k)$ of a space $(X,\Cal X)$ is called a canonical
kernel function (with respect to the probability measure $\mu$
on the space $(X,\Cal X)$) if
$$
\int f(x_1,\dots,x_{j-1},u,x_{j+1},\dots,x_k)\mu(\,du)=0\quad
\text{for all } 1\le j\le k \text{ \ and \ } x_s\in X,  \; s\neq j.
$$
Let me also introduce the notion of canonical functions in a more
general case, because this notion appears later in Proposition~5.
We call a function $f(x_1,\dots,x_k)$ on the $k$-fold product
$(X_1\times\cdots\times X_k, \Cal X_1\times\cdots\times \Cal X_k,
\mu_1\times\cdots\times \mu_k)$ of $k$ not necessarily identical
probability spaces $(X_j,\Cal X_j,\mu_j)$, $1\le j\le k$, if
$$
\int f(x_1,\dots,x_{j-1},u,x_{j+1},\dots,x_k)\mu_j(\,du)=0\quad
\text{for all } 1\le j\le k \text{ \ and \ } x_s\in X_s,  \; s\neq j.
$$
 
A $U$-statistic with a canonical kernel function is called degenerate.
Now I formulate Proposition~3.
\medskip\noindent
{\bf Proposition 3.} {\it Let us have a probability measure $\mu$ on
a space $(X,\Cal X)$ together with a sequence of independent and $\mu$
distributed random variables $\xi_1,\dots,\xi_n$ and an $L_2$-dense
class $\Cal F$ of canonical kernel functions $f=f(x_1,\dots,x_k)$
(with respect to the measure~$\mu$) with some parameter $D$ and
exponent $L$ on the product space $(X^k,\Cal X^k)$ which consists of
at most countably many functions, and satisfies conditions (1.2),
(1.3) and (1.4) with some $\sigma>0$. Let
$n\sigma^2>K((L+\beta)\log n+1)$ with a sufficiently large constant
$K=K(k)$. Then there exist some numbers $C=C(k)>0$,
$\gamma=\gamma(k)>0$ and threshold index $A_0=A_0(k)>0$ depending
only on the order $k$ of the $U$-statistics we consider such that the
degenerate $U$-statistics $I_{n,k}(f)$, $f\in\Cal F$, defined in (2.5)
satisfy the inequality
$$
P\(\sup_{f\in\Cal F}|n^{-k/2}I_{n,k}(f)|\ge A n^{k/2}\sigma^{k+1}\)
\le C e^{-\gamma A^{1/2k}n\sigma^2}\quad \text{if } A\ge A_0.
$$
} \medskip
(The constants in Propositions~2 and~3 may be different.) Before
deducing Proposition~2 from Proposition~3 I formulate a simple lemma
which will be useful also in the subsequent part of the paper. To
formulate it let us introduce the following notations.
 
Let some measure spaces $(Y_1,\Cal Y_1)$, $(Y_2,\Cal Y_2)$ and
$(Z,\Cal Z)$ be given together with a probability measure $\mu$
on the space $(Z,\Cal Z)$. Consider a function $f(y_1,z,y_2)$ on
the product space $(Y_1\times Z\times Y_2,\Cal Y_1\times\Cal Z
\times \Cal Y_2)$, $y_1\in\Cal Y_1$, $z\in\Cal Z$, $y_2\in\Cal Y_2$,
and define their projection
$$
P_\mu f(y_1,y_2)=\int f(y_1,z,y_2)\mu(\,dz),\quad y_1\in Y_1,\;
y_2\in Y_2, \tag2.6
$$
$$
\bar P_\mu f(y_1,z,y_2)=P_\mu f(y_1,y_2),
\quad y_1\in Y_1,\;z\in Z,\;y_2\in Y_2, \tag$2.6'$
$$
and
$$
\aligned
Q_\mu f(y_1,z,y_2)&=(I-\bar P_\mu) f(y_1,z,y_2)\\
&=f(y_1,z,y_2)-\bar P_\mu f(y_1,z,y_2), \quad
y_1\in Y_1,\;z\in Z,\;y_2\in Y_2.
\endaligned \tag$2.6''$
$$
(The difference between the operators $P_\mu$ and $\bar P_\mu$ is
that in the definition of the function $\bar P_\mu f$ we introduced
a fictive argument~$z$, i.e. $P_\mu f$ is defined on the space
$Y_1\times Y_2$ and $\bar P_\mu$ on the space $Y_1\times Z\times Y_2$.)
\medskip\noindent
{\bf Lemma~1.} {\it Let us have some measure spaces $(Y_1,\Cal Y_1)$,
$(Y_2,\Cal Y_2)$ and $(Z,\Cal Z)$, a probability measure $\mu$ on
the space $(Z,\Cal Z)$ and a probability measure $\rho$ on the product
space $(Y_1\times Y_2,\Cal Y_1\times\Cal Y_2)$. The transformations
$P_\mu$, $\bar P_\mu$ and $Q_\mu$ defined in (2.6)---($2.6''$) are
contractions from the space $L_2(Y_1\times Z\times Y_2,\rho\times\mu)$
to the spaces $L_2(Y_1\times Y_2,\rho)$ and
$L_2(Y_1\times Z\times Y_2,\rho\times\mu)$ respectively, i.e.
$$
\aligned
\|P_\mu f\|_{L_2,\rho}^2&=\int P_\mu
f(y_1,z,y_2)^2\rho(\,dy_1,\,dy_2) \\
&=\|\bar P_\mu f\|_{L_2,\rho\times\mu}^2
=\int\bar P_\mu f(y_1,z,y_2)^2\rho(\,dy_1,\,dy_2)\mu(\,dz)\\
&\le\|f\|_{L_2,\rho\times \mu}^2= \int
f(y_1,z,y_2)^2\rho(\,dy_1,\,dy_2)\mu(\,dz),
\endaligned \tag2.7
$$
and
$$
\aligned
\|Q_\mu f\|_{L_2,\rho}^2&=\int Q_\mu
f(y_1,z,y_2)^2\rho(\,dy_1,\,dy_2) \\
&=\int\(f(y_1,z,y_2)-\bar P_\mu
f(y_1,z,y_2)\)^2\rho(\,dy_1,\,dy_2)\mu(\,dz)\\
&\le\|f\|_{L_2,\rho\times \mu}^2= \int
f(y_1,z,y_2)^2\rho(\,dy_1,\,dy_2)\mu(\,dz).
\endaligned \tag$2.7'$
$$
 
If $\Cal F$ is an $L_2$-dense class of functions
$f(y_1,z,y_2)$ on the product space $(Y_1\times Z\times Y_2,\Cal
Y_1\times\Cal Z\times Y_2)$, $y_1\in\Cal Y_1$, $z\in\Cal Z$,
$y_2\in\Cal Y_2$ with parameter $D$ and
exponent $L$, then also the classes $\Cal F_\mu=\{P_\mu f,\:f\in \Cal
F\}$ and $\bar{\Cal F}_\mu=\{\bar P_\mu f,\:f\in \Cal F\}$ with the
functions $P_\mu f$ and $\bar P_\mu f$ defined in formulas (2.6) and
($2.6'$) are $L_2$-dense classes with parameter $D$ and exponent
$L$ in the spaces $(Y_1\times Y_2,\Cal Y_1\times\Cal Y_2)$ and
$(Y_1\times Z\times Y_2,\Cal Y_1\times\Cal Z\times\Cal Y_2$)
respectively, and the space $\Cal G_\mu=\{f-\bar P_\mu f,\;f\in
\Cal F\}$ defined in ($2.6''$) is an $L_2$-dense class with parameter
$2^LD$ and exponent~$L$ in the space $(Y_1\times Z\times Y_2,\Cal
Y_1\times\Cal Z\times\Cal Y_2$). Moreover, the class of functions
$\Cal G'_\mu=\{\frac12(f-\bar P_\mu f),\; f\in\Cal F\}$ is an
$L_2$-dense class with exponent $L$ and parameter~$D$.}
\medskip\noindent
{\it Proof of Lemma 1.}\/ The Schwarz inequality yields that
$P_\mu(f)^2\le\int f(y_1,z,y_2)^2\mu(\,dz)$, and the inequality
$\int [f(y_1,z,y_2)-\bar P_\mu f(y_1,z,y_2)]^2\mu(dz)\le \int
f(y_1,z,y_2)^2\mu(\,dz)$ also holds. Integrating these inequalities
with respect to the probability measure $\rho(\,dy_1,\,dy_2)$ we get
formulas (2.7) and ($2.7'$).
 
Let us consider an arbitrary probability measure $\rho$ on the space
$(Y_1\times Y_2,\Cal Y_1\times\Cal Y_2)$. To prove that $\Cal F_\mu$
is an $L_2$-dense class we have to find $m\le D \e^L$ functions
$f_j\in \Cal F_\mu$, $1\le j\le m$, such that
$\inff_{1\le j\le m}\int (f_j-f)^2\,d\rho\le \e^2$ for all $f\in
\Cal F_\mu$. But a similar property holds in the space $Y_1\times
Z\times Y_2$ with the probability measure $\rho\times\mu$. This
property together with the $L_2$ contraction property of $P_\mu$
formulated in (2.7) imply that $\Cal F_\mu$ is an  $L_2$-dense
class. The analogous property for $\bar{\Cal F}_\mu$ follows from
the already proved $L_2$-density property of $\Cal F_\mu$ and the
fact that by replacing a measure $\rho$ on $Y_1\times Z\times Y_2$
by the measure $\bar\rho\times\mu$, where $\bar\rho$ is the
projection of the measure $\rho$ to the space $Y_1\times Y$, i.e.\
$\bar\rho(B)=\rho(B\times Z)$ for $B\in \Cal Y_1\times\Cal Y_2$ we do
not change the $L_2$ norm of a difference $\bar P_\mu f-\bar P_\mu g$,
$f,g\in \Cal F$. Moreover, it equals to the $L_2$ norm of the
difference $P_\mu f-P_\mu g$ with respect to the measure $\bar\rho$.
Finally, the desired $L_2$-density property of the set $\Cal G_\mu$
can be deduced from the following observation. For any probability
measure $\rho$ on the space $Y_1\times Z\times Y_2$ and pair of
functions $f$ and $g$ such that $\int
(f-g)^2\frac12\(\,d\rho+\,d\bar\rho\times\,du\)\le \frac{\e^2}4$,
where $\bar\rho$ is the projection of the measure $\rho$ to the
space $Y_1\times Y_2$, $\int ((f-\bar P_\mu f)-(g-\bar P_\mu
g))^2\,d\rho\le 2\int (f-g)^2\,d\rho+2\int (\bar P_\mu f-\bar P_\mu
g)^2\,d\rho\le 2\int (f-g)^2\,d\rho+2\int
(f-g)^2\,d\bar\rho\times d\mu\le\e^2$.
This means that if $\{f_1,\dots,f_m\}$ is an $\frac\e2$-dense subset
of $\Cal F$ in the space $L_2(Y_1\times Z\times
Y_2,\Cal Y_1\times \Cal Z\times \Cal Y_2,\bar{\bar \rho})$ with
$\bar{\bar\rho}=\frac12(\rho+\bar\rho\times\mu)$, then
$\{Q_\mu f_1,\dots,Q_\mu f_m\}$ is an $\e$-dense subset of $\Cal
G_\mu$ in the space $L_2(Y_1\times Z\times
Y_2,\Cal Y_1\times \Cal Z\times \Cal Y_2,\rho)$. Moreover, if
$\{f_1,\dots,f_m\}$ is an $\e$-dense subset with respect to the
measure $\frac12\(\rho+\bar\rho\times\mu\)$, then
$\{\frac12Q_\mu f_1,\dots,\frac12Q_\mu f_m\}$ is an $\e$-dense subset
of $\Cal G_\mu'$ in the space $L_2(Y_1\times Z\times
Y_2,\Cal Y_1\times \Cal Z\times \Cal Y_2,\rho)$.
\medskip
 
To deduce Proposition~2 from Proposition~3 let us first introduce the
(random) probability measures $\mu^{(j)}$, $1\le j\le n$, concentrated
in the sample points $\xi_j$, i.e. let $\mu^{(j)}(A)=1$ if $\xi_j\in A$,
and $\mu^{(j)}(A)=0$ if $\xi_j\notin A$, $A\in \Cal A$. Then we can
write $\mu_n-\mu=\frac1n\(\summ_{j=1}^n\(\mu^{(j)}-\mu\)\)$, and
formula (1.1) can be rewritten as
$$
\align
J_{n,k}(f)=\dfrac1{n^{k/2}k!}\sum_{l=1}^k \sum_{j_l=1}^n  \int'
&f(u_1,\dots,u_k) \\
&\qquad \(\mu^{(j_1)}(\,du_1)-\mu(\,du_1)\)\dots
\(\mu^{(j_k)}(\,du_1)-\mu(\,du_1)\).
\endalign
$$
To rearrange the above sum in a way more appropriate for us let us
introduce the following notations: Let $\Cal P=\Cal P_k$ denote the
set of partitions of the set $\{1,2,\dots,k\}$, and given a sequence
$(j_1,\dots,j_k)$, $1\le j_s\le n$, $1\le s\le k$, of length $k$ let
$H(j_1,\dots,j_k)$ denote that partition of $\Cal P_k$ in which two
points $s$ and $t$, $1\le s,t\le k$, belong the same element of the
partition if $j_s=j_t$. Given a set $A$, let $|A|$ denote its
cardinality.
 
Let us rewrite the above expression for $J_{n,k}(f)$ in the form
$$
\align
J_{n,k}(f)=\dfrac1{n^{k/2}k!}\sum_{P\in \Cal P} &\sum\Sb (j_1,\dots,
j_k),\\ 1\le j_l\le n,\, 1\le l\le k \\ H(j_1,\dots,j_k)=P \endSb \int'
f(u_1,\dots,u_k)   \tag2.8 \\
&\qquad \(\mu^{(j_1)}(\,du_1)-\mu(\,du_1)\)\dots
\(\mu^{(j_k)}(\,du_1)-\mu(\,du_1)\).
\endalign
$$
 
Let us remember that the diagonals $u_s=u_t$, $s\neq t$, were
omitted from the domain of integration in the formula defining
$J_{n,k}(f)$. This implies that in the case $j_s=j_t$ the measure
$\mu^{(j_s)}(\,du_s)\mu^{(j_t)}(\,du_t)$ has zero measure in the
domain of integration. We have to understand the cancellation
effects caused by this relation. I want to show that because of
these cancellations the expression in formula (2.8) can be rewritten
as a linear combination of degenerated $U$-statistics with not too
large coefficients. The $U$-statistics taking part in this linear
combination can be bounded by means of Proposition~3, and this
yields an estimate sufficient for our purposes.  This seems to be a
natural approach, but the detailed proof demands some rather
unpleasant calculations.
 
Let us fix some $P\in\Cal P$ and investigate the inner sum at the
right-hand side of~(2.8) corresponding to this partition~$P$. For
the sake of simplicity let us first consider such a sum that
corresponds to a partition $P\in\Cal P$ which contains a set of the
form $\{1,\dots,s\}$ with some $s\ge2$. The products of measures
corresponding to the terms in the sum determined by such a
partition contain a part of length $s$ which has the form
$\(\mu^{(j)}(du_1)-\mu(du_1)\)\dots \(\mu^{(j)}(du_s)-\mu(du_s)\)$
with some $1\le j\le n$. This part of the product can be rewritten
in the domain of integration as
$$
\align
\summ_{l=1}^s &(-1)^{s-1}\mu(\,du_1)\dots\mu(\,du_{l-1})
(\mu^{(j)}(\,du_l)-\mu(\,du_l))\mu(\,du_{l+1})\dots\mu(\,du_s) \\
&\qquad +(-1)^{s-1}(s-1)\mu(du_1)\dots\mu(du_s).
\endalign
$$
Here we exploit that all other terms of this product disappears
in the domain of integration. Let us also observe that the term
$(-1)^{s-1}(s-1)\mu(du_1)\dots\mu(du_l)$ appears $n$-times as we
sum up for $1\le j\le n$. Similar calculation can be made for all
partitions $P\in\Cal P$ and all sets contained in the partitions,
only the notation of the indices will be more complicated. By
carrying out such a calculation the quantity $J_{n,k}(f)$ can be
rewritten as the linear combination of integrals of the function
$f(u_1,\dots,u_k)$ with respect to some product measures. The
components of these products of measures have either the form
$\(\mu^{(j_s)}(\,du_s)-\mu(\,du_s)\)$ or the form $\mu(\,du_s)$,
and all indices $j_s$ in a product are different. Let us observe
that to integrate a function $f$ with respect to
$\(\mu^{(j_s)}(\,du_s)-\mu(\,du_s)\)$ is the same as to apply
the operator $Q_{\mu,s}=I-\bar P_{\mu,s}$ for it and then to put
$\xi_s=u_s$ in the $s$-th argument of the function $Q_{\mu,s}f$,
and to integrate a function $f$ with respect to $\mu(du_s)$ is the
same as to apply the operator $P_{\mu,s}$ for it. Here $Q_{\mu,s}$
and $P_{\mu,s}$ are the  operators $Q_\mu$ and $P_\mu$ defined in
formulas ($2.6''$) and (2.6) if we choose in these formulas $Y_1$
as the product of the first $s-1$ components, $Z$ as the $s$-th
component and $Y_2$ as the product of the last $k-s$ components of
the $k$-fold product $X^k$.
 
Let us work out the details of the above indicated calculations and
for all sets $V\subset \{1,\dots,k\}$ let us gather in an internal
sum depending on $V$ those integrals for which the product of the
measures contain a component of the form
$\mu^{(j_s)}(\,du_s)-\mu(\,du_s))$, $1\le j_s\le n$,  if $s\in V$ and
a term $\mu(\,du_s)$ if $s\notin V$. In such a way we get the identity
$$
J_{n,k}(f)=\sum_{V\subset \{1,2,\dots,k\}} C(n,k,|V|)  n^{-|V|/2}
\frac1{k!}\sum\Sb 1\le j_s\le n,\\ J_s\neq j_{s'} \text{ if }s\neq s'
\\ \text{for } s\in V\endSb f_V(\xi_{j_s},s\in V) \tag2.9
$$
with the functions
$$
f_V(u_s,\,s\in V)=\prod_{s\in V}
Q_{\mu,s}  \!\!\!  \prod_{t\in\{1,\dots,k\}\setminus V}  \!\!\!
P_{\mu,t} \; f(u_1,u_2,\dots,u_k)\quad \text{for all }
V\subset\{1,\dots,k\} \tag2.10
$$
and some coefficients $C(n,k,|V|)$ which satisfy the inequality
$|C(n,k,|V|)|\le G(k)$ with some constant $G(k)>0$. The explicit
formula for $C(n,k,|V|)$ is rather complicated, but the above
estimate about the magnitude of this coefficient is sufficient for
our purposes. This estimate of $C(k,n,|V|)$ is sharp, because those
partitions $P\in\Cal P$  which contain the $|V|$ one-point subsets
of a set $V$ and $(k-|V|)/2$ subsets of cardinality 2 of
$\{1,\dots,k\}\setminus V$ yield a contribution of order
$n^{-k/2}n^{k/2-|V|/2}$ to the coefficient $C(n,k,|V|)n^{-|V|/2}$.
 
Let us observe that the inner sum corresponding to a set $V$ at the
right-hand side of (2.9) is a $U$-statistic with the kernel function
$f_V$ defined in (2.10). Hence to carry out our program we have to
understand the properties of this function $f_V$. It follows from
Lemma~1 that under the conditions of Proposition~1 the set of
functions $f_V$, $f\in \Cal F$, is an $L_2$-dense class with exponent
$L$ and parameter $2^{kL}D$, and
$\left|\int f^2_V(u_s,\,s\in V)\prodd_{s\in V}\mu(\,du_s)\right|\le
\sigma^2$ for all $V\in\{1,\dots,k\}$. Let me remark that this
estimate states in particular that the constant term $f_\emptyset$
defined in (2.10) with the choice $V=\emptyset$ satisfies the
inequality $|f_{\emptyset}|\le \sigma$. This estimate follows
directly from the Schwarz inequality, because
$$
f_{\emptyset}^2=\(\int f(u_1,\dots,u_k)\mu(\,du_1)\dots\mu(\,du_k)\)^2
\le\int f^2(u_1,\dots,u_k)\mu(\,du_1)\dots\mu(\,du_k)\le\sigma^2.
$$
 
Another important observation is that the functions $f_V$ are
canonical kernel functions with respect to the measure $\mu$. To
prove this statement let us observe that the canonical property of
a kernel function $f_V$ can be reformulated as
$P_{\mu,s}f_V(u_s,\,s\in V)=0$ for all $s\in V$ and sets of
parameters $u_t\in X$, $t\in V\setminus \{s\}$. This relation
follows from the observation that the operators $\bar P_{\mu,u}$,
$1\le u\le k$ are exchangeable, and $\bar P_{\mu,s}^2=\bar P_{\mu,s}$
which implies that $\bar P_{\mu,s}Q_{\mu,s}=\bar P_{\mu,s}(I-\bar
P_{\mu,s})=0$. (Actually, here we adapted the proof of the Hoeffding
decomposition of $U$-statistics to our case.)
 
Formula (2.9) yields that
$$
J_{n,k}(f)=\sum_{V\subset \{1,2,\dots,n\}} C(n,k,|V|)  n^{-|V|/2}
I_{n,|V|} (f_V(\xi_{j_s},s\in V)),
$$
and
$$
\align
&P\(\sup_{f\in\Cal F}|J_{n,k}(f)|\ge A n^{k/2}\sigma^{k+1}\) \\
&\qquad\qquad \le \sum_{V\subset \{1,2,\dots,n\}}
P\(\sup_{f\in\Cal F}| n^{-|V|/2}| I_{n,|V|}(f_V)|\ge \frac AT
n^{k/2}\sigma^{k+1}\)
\endalign
$$
with some appropriate constant $T=T(k)$. Observe that under the
Conditions of Proposition~2 $n\sigma^2\ge1$, hence
$n^{k/2}\sigma^{k+1}\ge n^{|V|/2}\sigma^{|V|+1}$. This means that if
the parameters $A_0$ and $K$ are sufficiently large in the conditions
of Propositions~2, then this conditions allow the application of
Proposition~3 to bound the probability
$P(n^{-|V|/2}|I_{n,|V|}(f_V)\ge \frac AT n^{k/2}\sigma^{k+1})
\le P(n^{-|V|/2}|I_{n,|V|}(f_V)\ge \frac AT n^{|V|/2}\sigma^{|V|+1})$
for all functions $f_V$. Thus we get that the inequality
$$
P\(\sup_{f\in\Cal F}|J_{n,k}(f)|\ge A n^{k/2}\sigma^{k+1}\)
\le C2^K e^{-\gamma(A/T)^{1/2k}n\sigma^2}
\le e^{-\gamma(A/2T)^{1/2k}n\sigma^2}
$$
holds for $A\ge A_0$ with some $T=T(k)$ if first the constant $K$
and then the constant $A_0$ are chosen sufficiently large in the
conditions of Proposition~2. This means that Proposition~3 implies
Proposition~2.
 
\beginsection
3. Some basic tools of the proof
 
First I formulate three results we apply in the proof of Proposition~3.
The first of them helps us to carry out some symmetrization arguments,
the second one yields a good estimate for the distribution of a
homogeneous polynomial  of independent random variables which
take values $\pm1$ with probability $\frac12$. Finally, the third
result enables us to reduce Proposition~3 to a simpler statement.
 
The first result, formulated in Lemma~2 is a slight generalization
of a simple lemma which can be found for instance in Pollard's
book~[9] ($8^\circ$ Symmetrization Lemma). I made this
generalization, because it is more appropriate for our purposes.
\medskip\noindent
{\bf Lemma 2. (Symmetrization Lemma)} {\it Let $Z(n)$ and $\bar Z(n)$,
$n=1,2,\dots$, be two sequences of random variables on a probability
space $(\Omega,\Cal A,P)$. Let a $\sigma$-algebra $\Cal B\subset
\Cal A$ be given on the probability space $(\Omega,\Cal A,P)$ together
with a $\Cal B$ measurable set $B$ and two numbers $\alpha>0$ and
$\beta>0$ such that the random variables $Z_n$, $n=1,2,\dots$ are
$\Cal B$ measurable, and the inequality
$$
P(|\bar Z_n|\le\alpha|\Cal B)(\oo)\ge\beta\quad \text{for all }
n=1,2,\dots \text{ if } \oo\in B \tag3.1
$$
holds.
Then
$$
P\(\sup_{1\le n<\infty}|Z_n|>\alpha+x\)\le\frac1\beta P\(\supp_{1\le
n<\infty}|Z_n-\bar Z_n|>x\)+(1-P(B))\quad\text{for all } x>0.
\tag3.2
$$
In particular, if the sequences $Z_n$, $n=1,2,\dots$, and $\bar Z_n$,
$n=1,2,\dots$, are two independent sequences of random  variables,
and $P(|Z_n|\le\alpha)\ge\beta$ for all $n=1,2,\dots$, then
$$
P\(\sup_{1\le n<\infty}|Z_n|>\alpha+x\)\le\frac1\beta P\(\supp_{1\le
n<\infty}|Z_n-\bar Z_n|>x\). \tag$3.2'$
$$
}\medskip\noindent
{\it Proof of Lemma 2.}\/ Put $\tau=\min\{n\: |Z_n|>\alpha+x)$ if
there exists such an $n$, and $\tau=0$ otherwise. Then
$$
\align
P(\{\tau=n\}\cap B)&\le \frac1\beta\int_{\{\tau=n\}\cap B} P(|\bar
Z_n|\le \alpha|\Cal B)\,dP
=\frac1\beta P(\{\tau=n\}\cap\{|\bar Z_n|\le\alpha\}\cap B)\\
&\le \frac1\beta P(\{\tau=n\}\cap\{|Z_n-\bar Z_n|>x\})
\quad \text{for all } n=1,2,\dots.
\endalign
$$
Hence
$$
\align
&P\(\sup_{1\le n<\infty}|Z_n|>\alpha+x\)-(1-P(B))\le
P\(\left\{\sup_{1\le n<\infty}|Z_n|>\alpha+x\right\}\cap B\) \\
&\qquad=\sum_{n=1}^\infty P(\{\tau=n\}\cap B)
\le \frac1\beta \sum_{n=1}^\infty P(\{\tau=n\}\cap\{|Z_n-\bar
Z_n|>x\}) \\
&\qquad \le\frac1\beta P\(\supp_{1\le n<\infty}|Z_n-\bar Z_n|>x\).
\endalign
$$
Thus formula 3.2 is proved. If $Z_n$ and $\bar Z_n$ are two
independent sequences, and $P(|Z_n|\le\alpha)\ge\beta$ for all
$n=1,2,\dots$, and we define $\Cal B$ as the $\sigma$-algebra
generated by the random variables $Z_n$, $n=1,2,\dots$, then the
condition (3.1) is satisfied also with $B=\Omega$. Hence relation
($3.2'$) holds in this case.  Lemma~2 is proved.
\medskip
The second result we need is a multi-dimensional version of
Hoeffding's inequality formulated in Proposition~A:
\medskip\noindent
{\bf Proposition A.} {\it  Let $\e_1,\dots,\e_n$  be independent
random variables, $P(\e_j=1)=P(\e_j=-1)=\dfrac12$, $1\le j\le n$.
Fix a positive integer~$k$ and define the random variable
$$
Z=\sum\Sb (j_1,\dots, j_k)\: 1\le j_l\le n \text{ for all } 1\le l\le
k\\ j_l\neq j_{l'} \text{ if }l\neq l' \endSb a(j_1,\dots, j_k)
\e_{j_1}\cdots \e_{j_k} \tag3.3
$$
with the help of some real numbers $a(j_1,\dots,j_k)$ which are given
for all sets of indices such that $1\le j_l\le n$, $1\le l\le k$, and
$j_l\neq j_{l'}$ if $l\neq l'$. Put
$$
S^2=\sum\Sb (j_1,\dots, j_k)\: 1\le j_l\le n \text{ for all } 1\le l\le
k\\ j_l\neq j_{l'} \text{ if }l\neq l' \endSb a^2(j_1,\dots, j_k)
\tag3.4
$$
Then
$$
P(|Z|>x)\le C \exp\left\{-B\(\frac xS\)^{2/k}\right\} \quad\text{for
all }x\ge 0 \tag3.5
$$
with some constants $B>0$ and $C>0$ depending only on the parameter
$k$. Relation (3.5) holds for instance with the choice
$B=\frac k{2e(k!)^{1/k}}$ and $C=e^k$.}
\medskip
Proposition~A is a relatively simple consequence of a famous and
important result of the probability theory, the so-called
hypercontractive inequality for Rademacher functions~(see~e.g.~[3]
or~[6]). The hypercontractive inequality yields some moment
inequalities that imply Proposition~A. Nevertheless, I did not find
this result in the literature. Therefore I explain in the Appendix
how it follows from the hypercontractive inequality.
\medskip\noindent
{\it Remark:}\/ The parameter $B$ given in Proposition~A is not
sharp. This is because the moment estimates I could prove are not
sharp enough. They are sufficient to give the right order of the
term in the exponent at the right-hand side of (3.5) but do not
give the best possible constant~$B$ in this estimate.
\medskip\noindent
 
Finally I formulate a decoupling type result which enables us to
reduce Proposition~3 to a similar but simpler statement. This
result compares the distribution of $U$-statistics with the
distribution of such systems whose coordinates are chosen
independently from each other. To make a clear distinction between
this object and usual $U$-statistics I shall call it independent
$U$-statistics. It is defined in the following way:
 
\medskip\noindent
{\bf Definition of independent $U$-statistics.}
{\it
Let us have $k$ independent copies $\xi_{1,s}$,\dots,~$\xi_{n,s}$,
$1\le s\le k$, of a sequence of independent and identically
distributed random variables $\xi_1,\dots,\xi_n$ with
distribution~$\mu$ on a measurable space $(X,\Cal X)$ together
with a function $f=f(x_1,\dots,x_k)$ on the $k$-th power $(X^k,
\Cal X^k)$ of the space $(X,\Cal X)$. We define with their help the
independent $U$-statistic~$\bar I_{n,k}(f)$ by the formula
$$
\bar I_{n,k}(f)=\frac1{k!}\summ\Sb 1\le j_s\le n,\; s=1,\dots, k\\
j_s\neq j_{s'} \text{ if } s\neq s'\endSb
f\(\xi_{j_1,1},\dots,\xi_{j_k,k}\). \tag3.6
$$
}\medskip
 
The following Proposition~B holds.
\medskip\noindent
{\bf Proposition B.} {\it Let us consider a countable sequence
$f_l(x_1,\dots,x_k)$, $l=1,2,\dots$, of functions on the $k$-fold
product $(X^k,\Cal X^k)$ of some space $(X,\Cal X)$ together with
some probability measure $\mu$ on the space $(X,\Cal X)$. Given a
sequence of independent and identically distributed random variables
$\xi_1,\xi_2,\dots$ with distribution~$\mu$ on $(X,\Cal X)$ together
with $k$ independent copies $\xi_{1,s},\xi_{2,s},\dots$,
$1\le s\le k$, of it we can define the $U$-statistics $I_{n,k}(f_l)$
and independent $U$-statistics $\bar I_{n,k}(f_l)$ for all
$l=1,2,\dots$ and $n=1,2,\dots$. They satisfy the inequality
$$
P\(\sup_{1\le l<\infty} \left| I_{n,k}(f_l)\right|>x\)\le
AP\(\sup_{1\le l<\infty}\left|\bar I_{n,k}(f_l)\right|>\gamma x\)
\tag3.7
$$
for all $x\ge0$ with some constants $A=A(k)>0$ and $\gamma
=\gamma(k)>0$ depending only on the order $k$ of the $U$-statistics.}
\medskip
 
I shall deduce Proposition~B from the result of paper~[5] of de la
Pe\~na and Mont\-go\-mery--Smith. At first sight one would think that
this result is not sufficient for our purposes, since it compares the
distribution function of a single $U$-statistic with its independent
$U$-statistic counterpart, i.e. the supremum with respect to a class
of functions is missing there. But this result is proved for general
Banach space valued random variables. Therefore, as I show below, its
application for an appropriate $L_\infty$ space yields the desired
result.
\medskip\noindent
{\it The proof of Proposition B (with the help of paper~[5].)}\/
Let us apply the first part of Theorem 1 of~[5] in the Banach
space $\ell_\infty$ consisting of the infinite sequences
$x=(x_1,x_2,\dots)$ of real numbers with norm $\|x\|=\!\!\supp_{1\le
l<\infty} \!|x_l|$ for the kernel functions
$f_{j_1,\dots,j_k}(x_1,\dots,x_k)=\bar f(x_1,\dots,x_k)$,
$\bar f=(f_1,f_2,\dots)$, mapping the space $(X^k,\Cal
X^k)$ into the space $\ell_\infty$. (Here we do not exploit that in
the result of~[5] the kernel functions may depend on the indices
$(j_1,\dots,j_k)$.) Then the result in~[5] states that
$$ \allowdisplaybreaks
\align
&P\(\left\|\summ\Sb 1\le j_s\le n,\; s=1,\dots, k\\
j_s\neq j_{s'} \text{ if } s\neq s'\endSb
\bar f\(\xi_{j_1},\dots,\xi_{j_k}\)\right\|>x\) \\
&\qquad \le AP\(\left\| \summ\Sb 1\le j_s\le n,\; s=1,\dots, k\\
j_s\neq j_{s'} \text{ if } s\neq s'\endSb
\bar f\(\xi_{j_1,1},\dots,\xi_{j_k,k}\)\right\|>\gamma x\)
\endalign
$$
with some universal constants $A=A(k)>0$ and $\gamma=\gamma(k)>0$,
and this statement is equivalent to relation (3.7).
\medskip\noindent
{\it Remark:}\/ Actually it would be enough to prove Proposition~B
only for the supremum of finitely many $U$-statistics with kernel
functions $f_1,\dots,f_N$ and then letting $N\to\infty$. In such a way
we can avoid the work with infinite dimensional Banach spaces. Such
an approach makes the proof simpler, in particular because some
measure-theoretical problems arise if we are working with Banach
spaces, where not all continuous linear functionals are measurable.
Such a difficulty really occurs if we are working with $L_\infty(X)$
spaces with a set $X$ of large cardinality. If we want to apply the
result of~[5] in the space $\ell_\infty$, then we have to check that it
is applicable in this case.
 
\medskip
Now I formulate the following Proposition~$3'$.
\medskip\noindent
{\bf Proposition~$3'$.} {\it Let us have a probability measure
$\mu$ on a space $(X,\Cal X)$ together with $k$ independent copies
$\xi_{1,s},\dots,\xi_{n,s}$, $1\le s\le k$, of a sequence of
independent and $\mu$ distributed random variables $\xi_1,\dots,\xi_n$
and a countable $L_2$-dense class $\Cal F$ of canonical kernel
functions $f=f(x_1,\dots,x_k)$ (with respect to the measure~$\mu$)
with some parameter $D$ and exponent $L$ on the product space
$(X^k,\Cal X^k)$ which satisfies conditions (1.2), (1.3) and (1.4)
with some $\sigma>0$. Let $n\sigma^2>K((L+\beta)\log n+1)$ with a
sufficiently large  constant $K=K(k)$. Then there exists some
threshold index $A_0=A_0(k)>0$ such that the independent
$U$-statistics $\bar I_{n,k}(f)$, $f\in\Cal F$, defined in (3.6)
satisfy the inequality
$$
P\(\sup_{f\in\Cal F}|n^{-k/2}\bar I_{n,k}(f)|\ge A n^{k/2}\sigma^{k+1}\)
\le e^{-A^{1/2k}n\sigma^2}\quad \text{if } A\ge A_0.
$$
} \medskip
 
Proposition~$3'$ and Proposition~B imply Proposition~3.
The proof of Proposition~$3'$ applies some ideas of a paper of
Alexander~[1]. Let me briefly explain them.
 
Let us restrict our attention to the case~$k=1$. In this case a
probability of the form $P\(n^{-1/2}\supp_{f\in\Cal F}
\left|\summ_{j=1}^n f(\xi_j)\right|>x\)$ has to be estimated. By
taking an independent copy of the sequence $\xi_n$ (which disappears
at the end of the of the calculation) a symmetrization argument can
be applied which reduces the problem to the estimation of the
probability $P\(n^{-1/2}\supp_{f\in\Cal F}\left|
\summ_{j=1}^n \e_jf(\xi_j)\right|>\bar x\)$, where the random
variables $\e_j$, $P(\e_j=1)=P(\e_j=-1)=\frac12$, $j=1,\dots,n$, are
independent, and they are independent also of the random variables
$\xi_j$. Beside this, the number $\bar x$ is only slightly smaller
than the number~$x/2$. Let us bound the conditional probability of
the event we have just introduced if the values random variables
$\xi_j$ are prescribed in it. This conditional probability can be
bounded by means of the one-dimensional version of Proposition~A,
and the estimate we get in such a way is useful if the conditional
variance of the random variable we have to handle has a good upper
bound. Such a bound exists, and some calculation reduces the original
problem to the estimation of the probability
$P\(n^{-1/2}\supp_{f\in\Cal F'}\left| \summ_{j=1}^n
f(\xi_j)\right|>x^{1+\alpha}\)$ with some new nice class of functions
$\Cal F'$ and number $\alpha>0$. This problem is very similar to the
original one, but it is simpler, since the number $x$ is replaced by
a larger number $x^{1+\alpha}$ in it. By repeating this argument
successively, in finitely many steps we get to an
inequality that clearly holds.
 
The above sketched argument suggests a backward induction procedure
to prove Proposition~$3'$. To carry out such a program I shall prove
a result formulated in Proposition~4. To do this first I introduce
the following notion.
\medskip\noindent
{\bf Definition of good tail behaviour for a class of $U$-statistics.}
{\it Let us have some measurable space $(X,\Cal X)$ and a probability
measure $\mu$ on it. Let us consider some class $\Cal F$ of functions
$f(x_1,\dots,x_k)$ on the $k$-fold product $(X^k,\Cal X^k)$ of the
space $(X,\Cal X)$. Fix some positive integer~$n$ and
positive number $\sigma>0$, and take $k$ independent copies
$\xi_{1,s},\dots,\xi_{n,s}$, $1\le s\le k$, of a
sequence of independent $\mu$-distributed random variables
$\xi_1,\dots,\xi_n$. Let us introduce with the help of these random
variables the independent $U$-statistics $\bar I_{n,k}(f)$, $f\in\Cal
F$. Given some real number $T>0$ we say that the set of independent
$U$-statistics determined by the class of functions $\Cal F$ has a
good tail behaviour at level~$T$ if the inequality
$$
P\(\sup_{f\in\Cal F}|n^{-k/2}\bar I_{n,k}(f)|\ge A
n^{k/2}\sigma^{k+1}\) \le \exp\left\{-A^{1/2k}n\sigma^2 \right\}
\quad \text{for all } A\ge T.  \tag3.8
$$
holds.}
\medskip
Now I formulate Proposition 4.
\medskip\noindent
{\bf Proposition 4.} {\it Let us fix a positive integer~$n$, real
number $\sigma>0$ and a probability measure $\mu$ on a measurable
space $(X,\Cal X)$ together with a countable $L_2$-dense class
$\Cal F$ of canonical kernel functions $f=f(x_1,\dots,x_k)$ (with
respect to the measure~$\mu$) on the $k$-fold product space
$(X^k,\Cal X^k)$ which has exponent $L$ and parameter~$D$, and the
number $D$ satisfies condition (1.4). Let us also assume that all
functions $f\in \Cal F$ satisfy the conditions
$\supp_{x_j\in X, 1\le j\le k}|f(x_1,\dots,x_k)|\le 2^{-(k+1)}$,
$\int f^2(x_1,\dots,x_k)\mu(\,dx_1)\dots\mu(\,dx_k)\le \sigma^2$,
and $n\sigma^2>K((L+\beta)\log n+1)$ with a sufficiently large fixed
number $K=K(k)$. There exists some real number $A_0=A_0(k)>1$ such
that for all classes of functions $\Cal F$ which satisfy the
conditions of Proposition~4 the sets of $U$-statistics determined by
the functions~$f\in\Cal F$ have a good tail behaviour at level~$T$
for some $T\ge A_0$, provided that they have a good tail behaviour
at level~$T^{4/3}$.} \medskip
It is not difficult to deduce Proposition~$3'$ from Proposition~4.
Indeed, let us observe that the set of $U$-statistics determined
by a class of functions $\Cal F$ satisfying the conditions of
Proposition~4 has a good tail-behaviour at level $n^{k/2}$, since
the probability at the left-hand side of (3.8) equals zero for
$u\ge n^{k/2}$. Then we get from Proposition~4 by induction with
respect to the number $j$, that this set of $U$-statistics has a
good tail-behaviour also for $T=n^{-(4/3)^jk/2}$ for $j=1,2,\dots$
if $n^{-(4/3)^jk/2}\ge A_0$. This implies that if a class of functions
$\Cal F$ satisfies the conditions of Proposition~4, then the set of
$U$-statistics determined by this class of functions has a good
tail-behaviour at level $T=A_0^{4/3}$, i.e. at a level which depends
only on the order $k$ of the (independent) $U$-statistics.
This result implies Proposition~$3'$, only we have to apply it not
directly for the class of functions~$\Cal F$ appearing in
Proposition~$3'$, but these functions have to be multiplied by a
sufficiently small positive number depending only on~$k$.
 
Thus to complete the proof of the Theorem it is enough to prove
Proposition~4. I describe its proof in the special case $k=1$ in
the next section. This case is considered separately, because it
may help to understand the ideas of the proof in the general case.
 
The main difficulty in the proof of Proposition 4 is related to a
symmetrization procedure which is an essential part of the proof.
We want to apply some randomization with the help of a
symmetrization argument, and this requires a special justification.
This is not a difficult problem in the case $k=1$, where it is
enough to calculate the variance of a $U$-statistic, but it becomes
hard for~$k\ge2$. In this case we have to give a good estimate on
certain conditional variances of some (independent) $U$-statistics
with respect to some appropriate conditions. To overcome this
difficulty we formulate a result in Proposition~5 and prove
Propositions~4 and~5 simultaneously. Their proof follows the
following line. First Proposition~4 and Proposition~5 are proved
for $k=1$. Then, if Propositions~4 and~5 are already proven for
all $k'<k$, then first we prove Proposition~4 for $k$, and
then Proposition~5 for the same~$k$. Proposition~5 has a similar
structure to Proposition~4. Before its formulation I introduce
the following definition.
\medskip\noindent
{\bf Definition of good tail behaviour for a class of integrals of
$U$-statistics.} {\it Let us have a product space $(X^k\times
Y,\Cal X^k\times\Cal Y)$ with some product measure $\mu^k\times\rho$,
where $(X^k,\Cal X^k,\mu^k)$ is the $k$-fold product of some
probability space $(X,\Cal X,\mu)$, and $(Y,\Cal Y,\rho)$ is some
other probability space. Fix some positive integer~$n$ and positive
number $\sigma>0$, and consider some class $\Cal F$ of functions
$f(x_1,\dots,x_k,y)$ on the product space $(X^k\times Y,\Cal
X^k\times\Cal Y,\mu^k\times\rho)$. Take $k$ independent copies
$\xi_{1,s},\dots,\xi_{n,s}$, $1\le s\le k$, of a sequence of
independent, $\mu$-distributed random variables $\xi_1,\dots,\xi_n$.
For all $f\in\Cal F$ and $y\in Y$ let us define the independent
$U$-statistics $\bar I_{n,k}(f,y)$ by means of these random variables
$\xi_{1,s},\dots,\xi_{n,s}$, $1\le s\le k$, and formula~(3.6).
Define with the help of these $U$-statistics $\bar I_{n,k}(f,y)$ the
random integrals
$$
H_{n,k}(f)=\int \bar I_{n,k}(f,y)^2\rho(\,dy), \quad f\in\Cal F.
\tag3.9
$$
Choose some real number $T>0$. We say that the set of random
integrals $H_{n,k}(f)$, $f\in\Cal F$, have a good tail behaviour at
level $T$ if
$$
P\(\sup_{f\in\Cal F} n^{-k}H_{n,k}(f)\ge A^2 n^k\sigma^{2k+2}\)
\le \exp\left\{-A^{1/(2k+1)}n\sigma^2 \right\}
\quad \text{for } A\ge T.
$$
}
\medskip\noindent
{\bf Proposition 5.} {\it Fix some positive integer $n$ and real
number $\sigma>0$, and let us have a product space $(X^k\times
Y,\Cal X^k\times\Cal Y)$ with some product measure $\mu^k\times\rho$,
where $(X^k,\Cal X^k,\mu^k)$ is the $k$-fold product of some
probability space $(X,\Cal X,\mu)$, and $(Y,\Cal Y,\rho)$ is some
another probability space. Let us have a countable $L_2$-dense class
$\Cal F$ of canonical functions $f(x_1,\dots,x_k,y)$ on the product
space $(X^k\times Y,\Cal X^k\times\Cal Y,\mu^k\times\rho)$ with some
exponent $L$ and parameter $D$ which satisfies condition~(1.4). Let
us also assume that the functions $f\in \Cal F$ satisfy the conditions
$$
\supp_{x_j\in X, 1\le j\le k, y\in Y}|f(x_1,\dots,x_k,y)|\le
2^{-(k+1)}
$$
and
$$
\int f^2(x_1,\dots,x_k,y)\mu(\,dx_1)\dots\mu(\,dx_k)\rho(\,dy)\le
\sigma^2 \quad  \text{for all } f\in \Cal F.
$$
Let the inequality $n\sigma^2>K((L+\beta)\log n+1)$ hold with a
sufficiently large fixed number $K=K(k)$.
 
There exists some number $A_0=A_0(k)>1$ such that for all
classes of functions $\Cal F$ which satisfy the conditions of
Proposition~5 the random integrals $H_{n,k}(f)$, $f\in\Cal F$,
defined in (3.9) have a good tail behaviour at level~$T$, provided
that they have a good tail behaviour at level~ $T^{(2k+1)/2k}$.}
\medskip
 
Similarly to the argument formulated after Proposition~4 an
inductive procedure yields the following corollary of Proposition~5.
\medskip\noindent
{\bf Corollary of Proposition 5.} {\it If the class of functions
$\Cal F$ satisfies the conditions of Proposition~5, then there
exists a constant $\bar A_0=\bar A_0(k)>0$ depending only on $k$
such that the integrals $H_{n,k}(f)$ determined by the class of
functions $\Cal F$ have a good tail behaviour at level $\bar A_0$.}
 
\medskip\noindent
{\bf 4. The proof of Proposition 4 in the case $k=1$}
 
\medskip\noindent
In this section Proposition~4 is proved in the special case $k=1$.
In this case we have to show that
$$
P\(\frac1{\sqrt n}\supp_{f\in\Cal F}\left|\summ_{j=1}^n
f(\xi_j)\right| \ge A n^{1/2}\sigma^{2}\) \le e^{-A^{1/2}
n\sigma^2} \quad \text{if } A\ge T
\tag4.1
$$
if we know the same estimate for $A>T^{4/3}$ and all classes of
functions satisfying the conditions of Proposition~4. This statement
will be proved by means of the following symmetrization argument.
\medskip\noindent
{\bf Lemma 3.} {\it Let the class of functions $\Cal F$ satisfy the
conditions of Proposition~4 for $k=1$. Let $\e_1,\dots,\e_n$ be a
sequence of independent random variables,
$P(\e_j=1)=P(\e_j=-1)=\frac12$, independent also of the
$\mu$ distributed random variables $\xi_1,\dots,\xi_n$. Then
$$
\aligned
&P\(\frac1{\sqrt n}\supp_{f\in\Cal F}\left|\summ_{j=1}^n
f(\xi_j)\right| \ge A
n^{1/2}\sigma^{2}\) \\
&\qquad \le 4P\(\frac1{\sqrt n}\supp_{f\in\Cal F}\left|\summ_{j=1}^n
\e_jf(\xi_j)\right| \ge \frac A3
n^{1/2}\sigma^{2}\) \quad\text{if } A\ge T.
\endaligned \tag4.2
$$
}\medskip\noindent
{\it Proof of Lemma 3.}\/ Let us construct an independent copy
$\bar\xi_1,\dots,\bar\xi_n$ of the sequence $\xi_1,\dots,\xi_n$ in
such a way that all three sequences $\xi_1,\dots,\xi_n$, \
$\bar\xi_1,\dots,\bar\xi_n$ and $\e_1,\dots,\e_n$ are independent.
Define the random variables $Z_n(f)=\frac1{\sqrt n}\summ_{h=1}^n
f(\xi_j)$ and $\bar Z_n(f)=\frac1{\sqrt n}\summ_{h=1}^n f(\bar\xi_j)$
for all $f\in\Cal F$. I claim that
$$
P\(\sup_{f\in\Cal F}|Z_n(f)|> A\sqrt n\sigma^2\)\le
2P\(\sup_{f\in\Cal F}|Z_n(f)-\bar Z_n(f)|> \frac23 A\sqrt
n\sigma^2\). \tag4.3
$$
This relation follows from Lemma~2 (the symmetrization lemma)
applied for the countable sets $Z_n(f)$ and $\bar Z_n(f)$,
$f\in\Cal F$, with $x=\frac23 A\sqrt n\sigma^2$ and
$\alpha=\frac13 A\sqrt n\sigma^2$, since the fields $Z_n(f)$
and $\bar Z_n(f)$ are independent, and $P(|Z_n(f)|\le \alpha)>\frac12$
for all $f\in\Cal F$. Indeed, $E\bar Z_n(f)^2\le\sigma^2$, thus
Chebishev's inequality implies that $P(|\bar Z_n(f)|\le\sqrt2\sigma)
\ge\frac12$ for all $f\in\Cal F$. On the other hand, we have assumed
that $n\sigma^2\ge K$ with some sufficiently large constant $K>0$.
Hence $\sigma\le\frac1{\sqrt K}\sqrt n\sigma^2$, and
$\sqrt2\sigma\le \alpha=\frac13 A\sqrt n\sigma^2$ if the
constant $K$ is chosen sufficiently large.
 
Let us observe that the random field
$$
Z_n(f)-\bar Z_n(f)=\frac1{\sqrt n}\sum_{j=1}^n \(f(\xi_j)
-f(\bar\xi_j)\), \quad f\in \Cal F,  \tag4.4
$$
and its randomization
$$
\frac1{\sqrt n}\sum_{j=1}^n \e_j \(f(\xi_j)
-f(\bar\xi_j)\), \quad f\in \Cal F,  \tag$4.4'$
$$
have the same distribution. Indeed, even the conditional distribution
of ($4.4'$) under the condition that the values of the $\e_j$-s are
prescribed agrees with the distribution of (4.4) for all possible
values of the $\e_j$-s. This follows from the observation that the
distribution of the field (4.4) does not change if we exchange the
random variables $\xi_j$ and $\bar\xi_j$ for certain indices $j$,
and this corresponds to considering the conditional distribution of
the field in ($4.4'$) under the condition that $\e_j=-1$ for these
indices $j$, and $\e_j=1$ for the remaining ones.
 
The above relation together with formula (4.3) imply that
$$
\align
&P\(\frac1{\sqrt n}\supp_{f\in\Cal F}\left|\summ_{j=1}^n
f(\xi_j)\right|  \ge A n^{1/2}\sigma^{2}\)\\
&\qquad \le 2P\(\frac1{\sqrt n}\supp_{f\in\Cal F}\left|\summ_{j=1}^n
\e_j\[f(\xi_j)-\bar f(\xi_j)\]\right| \ge\frac23 A
n^{1/2}\sigma^{2}\) \\
&\qquad\le 2P\(\frac1{\sqrt n}\supp_{f\in\Cal F}\left|\summ_{j=1}^n
\e_jf(\xi_j)\right| \ge\frac A3 n^{1/2}\sigma^{2}\) \\
&\qquad\qquad+ 2P\(\frac1{\sqrt n}\supp_{f\in\Cal F}\left|
\summ_{j=1}^n \e_jf(\bar\xi_j)\right| \ge\frac
A3n^{1/2}\sigma^{2}\) \\
&\qquad=4P\(\frac1{\sqrt n}\supp_{f\in\Cal F}\left|\summ_{j=1}^n
\e_jf(\xi_j)\right| \ge\frac A3n^{1/2}\sigma^{2}\)
\endalign
$$
Lemma~3 is proved.
\medskip
 
To prove Proposition~4 for $k=1$  let us investigate the conditional
probability
$$
P(f,A|\xi_1,\dots,\xi_n)=
P\(\left.\frac1{\sqrt n}\left|\summ_{j=1}^n
\e_jf(\xi_j)\right| \ge \frac
A6\sqrt n\sigma^2\right|\xi_1,\dots,\xi_n\)
$$
for all functions $f\in\Cal F$, $A\ge T$ and values
$(\xi_1,\dots,\xi_n)$. By Proposition~A (with $k=1$) we can write
$$
P(f,A|\xi_1,\dots,\xi_n)\le C\exp\left\{-\frac{\frac B{36}
A^2 n\sigma^4}{S^2(f,\xi_1,\dots,\xi_n)}\right\} \tag4.5
$$
with
$$
S^2(f,x_1,\dots,x_n)=\frac1n\sum_{j=1}^n f^2(x_j), \quad f\in \Cal F.
$$
Let us introduce the set
$$
H=H(A)=\left\{(x_1,\dots,x_n)\: \sup_{f\in F} S^2(f,x_1,\dots,x_n)\ge
\(1+A^{4/3}\)\sigma^2\right\}. \tag4.6
$$
I claim that
$$
P((\xi_1,\dots,\xi_n)\in H)\le e^{-A^{2/3} n\sigma^2}\quad\text{ if }
A\ge T. \tag$4.6'$
$$
To prove relation ($4.6'$) let us consider the functions $\bar f=\bar
f(f)$ for all $f\in \Cal F$ defined by the formula $\bar
f(x)=f^2(x)-\int f^2(x)\mu(\,dx)$, and introduce the class of functions
$\Cal F'=\{\bar f(f)\: f\in\Cal F\}$. Let us show that the class of
functions $\Cal F'$ satisfies the conditions of Proposition~4, hence
the estimate (4.1) holds for the class of functions $\Cal F'$
if $A\ge T^{4/3}$.
 
The relation $\int \bar f(x)\mu(\,dx)=0$ clearly holds. (In the case
$k=1$ this means that $\bar f$ is a canonical function.) The condition
$\sup| \bar f(x)|\le\frac 18<\frac14$ also holds if $\sup |f(x)|\le
\frac14$, and $\int \bar f^2(x)\mu(\,dx)\le \int f^4(x)\mu(\,dx)\le
\frac 14\int f^2(x)\,\mu(\,dx)\le\frac{\sigma^2}4<\sigma^2$ if $f\in
\Cal F$. It remained to show that $\Cal F'$ is an $L_2$-dense class
with exponent $L$ and parameter $D$.
 
To show this observe that $\int (\bar f(x)-\bar g(x))^2\rho(\,dx)\le
2\int(f^2(x)-g^2(x))^2\rho(\,dx)+
2\int(f^2(x)-g^2(x))^2\mu(\,dx)\le2 (\supp (|f(x)|+|g(x)|)^2
\(\int (f(x)-g(x))^2(\rho(\,dx)+\mu(\,dx)\)\le  \int
(f(x)-g(x))^2\bar\rho(\,dx)$ for all $f, g\in\Cal F$, $\bar f=\bar
f(f)$, $\bar g=\bar g(g)$ and probability measure $\rho$, where
$\bar\rho=\frac{\rho+\mu)}2$. This means that if $\{f_1,\dots,f_m\}$
is an $\e$-dense subset of $\Cal F$ in the space $L_2(X,\Cal
X,\bar\rho)$, then $\{\bar f_1,\dots,\bar f_m\}$ is an $\e$-dense
subset of $\Cal F'$ in the space $L_2(X,\Cal X,\rho)$, and
not only $\Cal F$, but also $\Cal F'$ is an $L_2$-dense class with
exponent $L$ and parameter $D$.
 
We get, by applying formula (4.1) for the number $A^{4/3}\ge
T^{4/3}$ and the class of functions $\Cal F'$ that
$$
\align
P((\xi_1,\dots,\xi_n)\in H)&=P\(\sup_{f\in F} \(\frac1n \sum_{h=1}^n
\bar f(\xi_j) +\frac1n \sum_{h=1}^n E f^2(\xi_j)\)
\ge \(1+A^{4/3}\)\sigma^2\)\\
&\le P\(\sup_{f\in F} \frac1{\sqrt n} \sum_{h=1}^n \bar f(\xi_j)
\ge A^{4/3}n^{1/2}\sigma^2\) \le e^{-A^{2/3} n\sigma^2},
\endalign
$$
i.e. relation ($4.6'$) holds.
 
Formula (4.5) and the definition (4.6) of the set $H$ yield the estimate
$$
P(f,A|\xi_1,\dots,\xi_n)\le Ce^{- B A^{2/3} n\sigma^2/40} \quad
\text{if }(\xi_1,\dots,\xi_n)\notin H \tag4.7
$$
for all $f\in \Cal F$ and $A\ge T$ for the conditional
probability $P(f,A|\xi_1,\dots,\xi_n)$.
Let us introduce the conditional probability
$$
P(\Cal F,A|\xi_1,\dots,\xi_n)=
P\(\left.\sup_{f\in \Cal F} \frac1{\sqrt n}\left|\summ_{j=1}^n
\e_jf(\xi_j)\right| \ge \frac
A3\sqrt n\sigma^2\right|\xi_1,\dots,\xi_n\)
$$
for all $(\xi_1,\dots,\xi_n)$ and $A\ge T$. We shall
estimate this conditional probability with the help of relation (4.7)
if $(\xi_1,\dots,\xi_n) \notin H$. Given some set of $n$~points
$(x_1,\dots,x_n)$ in the space $(X,\Cal X)$ let us introduce the
measure $\nu=\nu(x_1,\dots,x_n)$ on $(X,\Cal X)$ in such a way that
$\nu$ is concentrated in the points $x_1,\dots,x_n$, and
$\nu(\{x_j\})=\frac1n$. If $\int f^2(u)\nu(\,du)\le\delta^2$ for a
function $f$, then $\left|\frac1{\sqrt n}\summ_{j=1}^n
\e_jf(x_j)\right|\le n^{1/2}\int|f(u)|\nu(\,du)\le n^{1/2}\delta$.
Since we have assumed that $n\sigma^2\ge1$, this estimate implies
that if $f$ and $g$ are two functions such that $\int
(f-g)^2\nu(\,dx)\le \delta^2$ with $\delta=\frac A{6n}$,
then $\left|\frac1{\sqrt n}\summ_{j=1}^n \e_jf(x_j)-
\frac1{\sqrt n}\summ_{j=1}^n \e_jg(x_j)\right|\le\frac A{6\sqrt n}
\le\frac A6 \sqrt n\sigma^2$.
 
Given some (random) point $(\xi_1,\dots,\xi_n)\in H$ let us consider
the measure $\nu=\nu(\xi_1,\dots,\xi_n)$ corresponding to it, and
choose a $\bar\delta$-dense subset $\{f_1,\dots,f_m\}$ of $\Cal F$ in
the space $L_2(X,\Cal X,\nu)$ with $\bar\delta=\frac1{6n}\le\delta=
\frac A{6n}$, whose cardinality $m$ satisfies the inequality $m\le
D\bar\delta^{-L}$. This is possible because of the $L_2$-dense
property of the class~$\Cal F$. (This is the point where the
$L_2$-dense property of the class of functions $\Cal F$ is exploited
in its full strength.) The above facts imply that $P(\Cal
F,A|\xi_1,\dots,\xi_n)\le\summ_{l=1}^m P(f_l,A|\xi_1,\dots,\xi_n)$ with
these functions $f_1,\dots,f_m$. Hence relation (4.7) yields that
$$
P(\Cal F,A|\xi_1,\dots,\xi_n)\le CD(6n)^Le^{- B
A^{2/3} n\sigma^2/40} \quad \text{if }(\xi_1,\dots,\xi_n)\notin H
\text{ and } A\ge T.
$$
This inequality together with Lemma~3 and estimate~($4.6'$) imply that
$$
\aligned
&P\(\frac1{\sqrt n}\supp_{f\in\Cal F}\left|\summ_{j=1}^n
f(\xi_j)\right| \ge A n^{1/2}\sigma^{2}\)
\le 4P\(\frac1{\sqrt n}\supp_{f\in\Cal F}\left|\summ_{j=1}^n
\e_jf(\xi_j)\right| \ge \frac A3 n^{1/2}\sigma^{2}\)\\
&\qquad \le 4CD(6n)^Le^{- B A^{2/3}n\sigma^2/40}
+4e^{-A^{2/3}n\sigma^2} \quad \text{if } A\ge T.
\endaligned \tag4.8
$$
Since we have a better power of $A$ in the exponent at the
right-hand side of formula (4.8) than we need, the relation
$n\sigma^2\ge K((L+\beta)\log n+1)$ holds, and we have the right to
choose the constants $K$ and $A_0$, $A\ge A_0$, sufficiently large,
it is not difficult to deduce relation (3.8) from relation (4.8).
Indeed, the expression in the exponent at the right-hand side of
(4.8) satisfies the inequality $\frac B{40} A^{2/3} n\sigma^2\ge
A^{1/2} n\sigma^2+K((L+\beta)\log n+1)$ if $A_0$ is sufficiently
large, and
$$
\align
P&\(\frac1{\sqrt n}\supp_{f\in\Cal F}\left|\summ_{j=1}^n
f(\xi_j)\right| \ge An^{1/2}\sigma^{2}\)\\
&\qquad \le 4C(6n)^{\beta+L}e^{-K} n^{-K(L+\beta)}
e^{-A^{1/2} n\sigma^2}+4e^{-A^{2/3}n\sigma^2}\le e^{-A^{1/2}n\sigma^2}
\endalign
$$
if $A\ge T$, and the constants $A_0$ and $K$ are chosen sufficiently
large.
 
\medskip\noindent
{\bf 5. The symmetrization argument}
 
\medskip\noindent
In the proof of Propositions~4 and~5 we need two symmetrization
results for all $k\ge1$ which play the same role as Lemma~3 in the
case $k=1$. These results are described in Lemmas~4A and~4B. In
this section these results are formulated and proved. The proofs
go by induction with respect to $k$. During the proof of Propositions~4
and~5 for~$k$ we may assume that they hold for $k'<k$.
\medskip\noindent
{\bf Lemma 4A.} {\it Let $\Cal F$ be a class of functions on the
space $(X^k,\Cal X^k)$ which satisfies the conditions of Proposition~4
with some probability measure $\mu$. Let us have $k$ independent
copies $\xi_{1,s},\dots,\xi_{n,s}$, $1\le s\le k$, of a sequence of
independent $\mu$ distributed random variables $\xi_1,\dots,\xi_n$,
and a sequence of independent random variables $\e=(\e_1,\dots,\e_n)$,
$P(\e_s=1)=P(\e_2=-1)=\frac12$, which is independent also of the
random variables $\xi_{j,s}$, $1\le j\le n$, $1\le s\le k$. Consider
the independent $U$-statistics $\bar I_{n,k}(f)$, $f\in\Cal F$,
defined from these random variables by formula~(3.6) and their
randomized version
$$
\bar I_{n,k}^{\e}(f)=\frac1{k!}\summ\Sb 1\le j_s\le n,\; s=1,\dots,
k\\ j_s\neq j_{s'} \text{ if } s\neq s'\endSb
\e_{j_1}\cdots\e_{j_k}f\(\xi_{j_1,1},\dots,
\xi_{j_k,k}\),  \quad f\in\Cal F. \tag5.1
$$
There exists some constant $A_0=A_0(k)$ such that the inequality
$$
\aligned
P\(\sup_{f\in\Cal F} n^{-k/2}\left|\bar
I_{n,k}(f)\right|>An^{k/2}\sigma^{k+1}\)&<
2^{k+1}P\(\sup_{f\in\Cal F} \left|\bar I_{n,k}^{\e}(f)\right|
>2^{-(k+1)}A n^k\sigma^{k+1}\)\\
&\qquad+2^kn^{k-1}e^{-A^{1/(2k-1)} n\sigma^2/k}
\endaligned \tag 5.2
$$
holds for all $A\ge A_0$.}
\medskip
Before formulating Lemma 4B needed in the proof of Proposition~5
I introduce some notations. Some of them will be needed later.
 
Let us consider a set of functions $\Cal F$ of functions
$f(x_1,\dots,x_k,y)\in \Cal F$ on a space $(X^k\times Y, \Cal X^k
\times \Cal Y,\mu^k\times\rho)$ which satisfies the conditions of
Proposition~5. Let us choose $2k$ independent copies
$\xi_{1,s}^{(1)},\dots,\xi_{n,s}^{(1)}$,
$\xi_{1,s}^{(-1)},\dots,\xi_{n,s}^{(-1)}$, $1\le s\le k$, of a
sequence of independent $\mu$ distributed random variables
$\xi_1,\dots,\xi_k$ together with a sequence of independent random
variables $(\e_1,\dots,\e_n)$, $P(e_s=1)=P(\e_s-1)=\frac12$, $1\le
s\le n$ which are independent of them. For all subsets
$V\subset\{1,\dots,k\}$
of the set $\{1,\dots,k\}$ let $|V|$ denote the cardinality of this set,
and define for all functions $f(x_1,\dots,x_k,y)\in \Cal F$ and
$V\subset\{1,\dots,k\}$ the independent $U$-statistics
$$
\bar I_{n,k}^V(f,y)=\frac1{k!}\summ\Sb 1\le j_s\le n,\; s=1,\dots, k\\
j_s\neq j_{s'} \text{ if } s\neq s'\endSb
f\(\xi_{j_1,1}^{(\delta_1)},\dots,\xi_{j_k,k}^{(\delta_k)},y\),\quad
f\in\Cal F, \tag5.3
$$
where $\delta_s=\pm1$, $1\le s\le k$, $\delta_s=1$ if $s\in V$,
and $\delta_s=-1$ if $s\notin V$, together with the random variables
$$
H_{n,k}^V(f)=\int \bar I_{n,k}^V(f,y)^2\rho(\,dy), \quad f\in\Cal
F. \tag$5.3'$
$$
Put
$$
\bar I_{n,k}(f,y)=\bar I_{n,k}^{\{1,\dots,k\}}(f,y),\quad
H_{n,k}(f)=H_{n,k}^{\{1,\dots,k\}}(f), \tag$5.3''$
$$
i.e. these random variables appear if $V=\{1,\dots,k\}$ is taken in
the previous definitions, and the random variables $\xi_{j,s}^{(1)}$,
$1\le j\le n$, $1\le s\le k$ are inserted in the formulas defining
these random variables.
 
Let us also define the `randomized version' of the random variables
$\bar I_{n,k}^V(f,y)$ and $H_{n,k}^V(f)$ as
$$
\bar I_{n,k}^{(V,\e)}(f,y)=\frac1{k!}\summ\Sb 1\le j_s\le n,\;
s=1,\dots, k\\ j_s\neq j_{s'} \text{ if } s\neq s'\endSb
\e_{j_1}\cdots\e_{j_k}f\(\xi_{j_1,1}^{(\delta_1)},\dots,
\xi_{j_k,k}^{(\delta_k)},y\),\quad f\in\Cal F, \tag5.4
$$
where $\delta_s=1$ if $s\in V$, and $\delta_s=-1$ if $s\notin V$, and
$$
H_{n,k}^{(V,\e)}(f)=\int \bar I_{n,k}^{(V,\e)}(f,y)^2\rho(\,dy)
,\quad f\in\Cal F. \tag$5.4'$
$$
 
Let us also introduce the random variables
$$
\bar W(f)=\int\[\sum_{V\subset \{1,\dots,k\}} (-1)^{|V|}\bar
I_{n,k}^{(V,\e)}(f,y)\]^2\rho(\,dy), \quad f\in\Cal F. \tag5.5
$$
Now I formulate the symmetrization result applied
in the proof of Proposition~5.
\medskip\noindent
{\bf Lemma 4B.} {\it Let $\Cal F$ be a set of functions on
$(X^k\times Y,\Cal X^k\times\Cal Y)$ which satisfies the conditions
of Proposition~5 with some probability measure $\mu^k\times\rho$.
Let us have $2k$ independent copies
$\xi_{1,s}^{\pm1},\dots,\xi_{n,s}^{\pm1}$, $1\le s\le k$, of a
sequence of independent $\mu$ distributed random variables
$\xi_1,\dots,\xi_n$ together with a sequence of independent random
variables $\e_1,\dots,\e_n$, $P(e_s=1)=P(\e_s=-1)=\frac12$, $1\le s\le
n$, which is independent also of the previously considered sequences.
 
There exists some $A_0=A_0(k)$ such that if the integrals
$H_{n,k}(f)$, $f\in\Cal F$, determined by this class of functions
$\Cal F$ have a good tail behaviour at level $T^{(2k+1)/2k}$ for
some $T\ge A_0$, (this property was defined at the end of
Section~3), then the inequality
$$
\aligned
P\(\sup_{f\in\Cal F} H_{n,k}(f)>A^2n^{2k}\sigma^{2(k+1)}\)
&<2P\(\sup_{f\in\Cal F} \left|\bar W(f)\right|
>\frac{A^2}2 n^{2k}\sigma^{2(k+1)}\)\\
&\qquad+2^{2k+1}n^{k-1}e^{-A^{1/2k} n\sigma^2/k}
\endaligned \tag 5.6
$$
holds with the random variables $H_{n,k}(f)$ and $\bar W(f)$ defined
in formulas $(5.3'')$ and (5.5) for all $A\ge T$.}
\medskip
 
Let us observe that in the symmetrization argument of Lemma~4B we have
applied the symmetrization $\bar I_{n,k}^{(V,\e)}(f,y)$ of  $\bar
I_{n,k}^(V(f,y)$, (compare formulas (5.3) and (5.4)), and compared
the integral of the square of the random function $\bar I_{n,k}(f,y)$
with the integral of the square of a linear combination of the random
functions $\bar I_{n,k}^{(V,\e)}(f,y)$. After this integration the
effect of the `randomizing factors' $\e_j$ will be weaker.
Nevertheless, also such an estimate will be sufficient for us. But
the effect of this symmetrization procedure has to be followed more
carefully. Hence a corollary of Lemma~4B will be presented which can
be better applied than the original lemma. We get it by rewriting the
random variable $\bar W(f)$ defined in (5.5) in another form with the
help of some diagrams introduced below.
 
Let $\Cal G=\Cal G(k)$ denote the set of all diagrams consisting of
two rows such that both rows are the set $\{1,\dots,k\}$ and the
diagrams of $\Cal G$ contain some edges $(l_1,l_1')$,\dots,
$(l_s,l_s')$, $0\le s\le k$ connecting some points (vertices) of
the first row with some point (vertex) of the second row. The
vertices $l_1,\dots,l_s$ in the first row are all different, and
the same relation holds also for the vertices $l_1',\dots,l_s'$ in
the second row. For each diagram $G\in\Cal G$ let us define
$e(G)=\{(l_1,l_1')\dots,(l_s,l_s')\}$, the set of its edges,
$v_1(G)=\{l_1,\dots,l_s\}$, the set of its vertices in the first
row and $v_2(G)=\{l_1',\dots,l_s'\}$, the set of its vertices in the
second row.
 
Given some diagram $G\in \Cal G$ and two sets
$V_1,V_2\subset\{1,\dots,k\}$, we define with the help of the
random variables $\xi_{s,1}^{(1)},\dots,\xi_{s,n}^{(1)}$,
$\xi_{s,1}^{(-1)},\dots,\xi_{s,n}^{(-1)}$, $1\le s\le k$, and
$\e=(\e_1,\dots,\e_n)$ taking part in the definition of the
expressions $\bar W(f)$, $f\in\Cal F$, the random variables
$H_{n,k}(f|G,V_1,V_2)$:
$$
\aligned
H_{n,k}(f|G,V_1,V_2)&=\sum\Sb(j_1,\dots,j_k,\;j'_1,\dots,j'_k) \\
1\le j_s\le n,\, j_s\neq j_{s'} \text{ if }s\neq s',\,1\le s,s'\le k,\\
1\le j'_s\le n,\, j'_s\neq j'_{s'}\text { if } s\neq s',\,1\le s,s'\le
k,\\ j_s=j'_{s'} \text { if } (s,s')\in e(G),\; j_s\neq j'_{s'} \text
{ if } (s,s')\notin e(G)\endSb \prod_{s\notin v_1(G)}\e_{j_s}
\prod_{s\notin v_2(G)}\e_{j'_s} \\
&\qquad \int
f(\xi_{j_1,1}^{(\delta_1)},\dots,\xi_{j_k,k}^{(\delta_k)},y)
f(\xi_{j'_1,1}^{(\bar\delta_1)},\dots,\xi_{j'_k,k}^{(\bar\delta_k)},y)
\rho(\,dy), \quad f\in\Cal F,
\endaligned \tag5.7
$$
where $\delta_s=1$ if $s\in V_1$, $\delta_s=-1$ if $s\notin V_1$,
and $\bar\delta_s=1$ if $s\in V_2$, $\bar\delta_s=-1$ if $s\notin V_2$.
 
With the help of these random variables we can write that
$$
\bar W(f)=\sum_{G\in \Cal G,\, V_1,V_2\subset \{1,\dots,k\}}
(-1)^{|V_1|+|V_2|} H_{n,k}(f|G,V_1,V_2) \quad \text{for all }
f\in\Cal F,
$$
because
$$
\int\bar I_{n,k}^{(V_1,\e)}(f,y)\bar I_{n,k}^{(V_2,\e)}(f,y)\rho(\,dy)
=\sum_{G\in\Cal G} H_{n,k}(f|G,V_1,V_2), \quad \text{for all }
V_1,V_2\subset\{1,\dots,k\}.
$$
 
Since the number of terms in this sum is less than $2^{4k}k!$, it
implies that Lemma~4B has the following corollary:
\medskip\noindent
{\bf Corollary of Lemma 4B.} {\it Let a set of functions $\Cal F$
satisfy the conditions of Proposition~5. Then there exists some
$A_0=A_0(k)$ such that if the integrals
$H_{n,k}(f)$, $f\in\Cal F$, determined by this class of functions
$\Cal F$ have a good tail behaviour at level $T^{(2k+1)/2k}$ for
some $T\ge A_0$, then the inequality
$$
\align
&P\(\sup_{f\in\Cal F} H_{n,k}(f)>A^2n^{2k}\sigma^{2(k+1)}\)\\
&\qquad\qquad\le 2\sum_{G\in \Cal G,\, V_1,V_2\in\{1,\dots,k\}}
P\(\sup_{f\in\Cal F} \left |H_{n,k}(f|G,V_1,V_2)\right|
>\frac{A^2}{2^{4k+1}k!} n^{2k}\sigma^{2(k+1)}\) \\
&\qquad\qquad\qquad +2^{2k+1}n^{k-1}e^{-A^{1/2k} n\sigma^2/k} \tag 5.8
\endalign
$$
holds with the random variables $H_{n,k}(f)$ and $H_{n,k}(f|G,V_1,V_2)$
defined in formulas $(5.3'')$ and (5.7) for all $A\ge T$.}
\medskip
 
The proof of Lemmas 4A and 4B uses the result of the following
Lemma~5 which states that certain random vectors have the same
distribution.
\medskip\noindent
{\bf Lemma 5.} {\it Let $\e=(\e_1,\dots,\e_n)$ be a sequence of
independent random variables, $P(\e_s=1)=P(\e_s=-1)=\frac12$,
$1\le s\le n$, which is independent also of $2k$ fixed independent
copies $\xi_{1,s}^{(1)},\dots,\xi_{n,s}^{(1)}$ and
$\xi_{1,s}^{(-1)},\dots,\xi_{n,s}^{(-1)}$, $1\le s\le k$, of a
sequence $\xi_1,\dots,\xi_n$ of independent $\mu$ distributed random
variables. \medskip
\item a) Let $\Cal F$ be a class of functions which satisfies the
conditions of Proposition 4. With the help of the above
random variables introduce the independent $U$-statistic
$$
\bar I_{n,k}^V(f)=\frac1{k!}\summ\Sb 1\le j_s\le n,\; s=1,\dots, k\\
j_s\neq j_{s'} \text{ if } s\neq s'\endSb
f\(\xi_{j_1,1}^{(\delta_1)},\dots,\xi_{j_k,k}^{(\delta_k)}\),\quad
f\in\Cal F, \tag5.9
$$
for all sets $V\subset\{1,\dots,k\}$ and functions $f\in \Cal F$
together with its `randomized version'
$$
\bar I_{n,k}^{(V,\e)}(f)=\frac1{k!}\summ\Sb 1\le j_s\le n,\;
s=1,\dots, k\\ j_s\neq j_{s'} \text{ if } s\neq s'\endSb
\e_{j_1}\cdots\e_{j_k}f\(\xi_{j_1,1}^{(\delta_1)},\dots,
\xi_{j_k,k}^{(\delta_k)}\),  \quad f\in\Cal F, \tag$5.9'$
$$
where $\delta_s=\pm1$, $1\le s\le k$, $\delta_s=1$ if $s\in V$, and
$\delta_s=-1$ if $s\notin V$.
 
Then the sets of random variables
$$
S(f)=\sum_{V\subset \{1,\dots,k\}} (-1)^{|V|}\bar I_{n,k}^V(f),
\quad f\in\Cal F, \tag5.10
$$
and sets of random variables
$$
\bar S(f)=\sum_{V\subset \{1,\dots,k\}} (-1)^{|V|}\bar
I_{n,k}^{(V,\e)}(f), \quad f\in\Cal F, \tag$5.10'$
$$
have the same joint distribution.
\medskip
\item b) Let $\Cal F$ be the class of functions satisfying
Proposition 5. For all functions $f\in \Cal F$ and
$V\subset\{1,\dots,k\}$ consider the independent $U$-statistics
determined by the random variables
$\xi_{1,s}^{(1)},\dots,\xi_{n,s}^{(1)}$ and
$\xi_{1,s}^{(-1)},\dots,\xi_{n,s}^{(-1)}$, $1\le s\le k$ by
formula (5.3), and define with their help the random variables
$$
W(f)=\int\[\sum_{V\subset \{1,\dots,k\}} (-1)^{|V|}\bar
I_{n,k}^V(f,y)\]^2\rho(\,dy), \quad f\in\Cal F. \tag5.11
$$
Then the random vectors $\{W(f)\: f\in \Cal F\}$ defined in (5.11)
and $\{\bar W(f)\: f\in \Cal F\}$ defined in (5.5)
have the same distribution.}
\medskip\noindent
{\it Proof of Lemma 5.} Let us consider Part a) of Lemma~5. I claim
that for all $M\in\{1,\dots,n\}$  the conditional distribution of
the random vector in $(5.10')$ under the condition that $\e_j=1$ if
$j\in M$ and $\e_j=-1$ if $\e_j\in\{1,\dots,n\}\setminus M$ agrees
with the distribution of the vector in (5.10). Since the
distribution of the vector in (5.10) does not change if we exchange
the random variables $\xi_{j,s}^{(1)}$ and $\xi_{j,s}^{(-1)}$ in it if
$j\notin M$, $1\le s\le k$, and do not exchange them otherwise, it is
enough to understand that the random vector we get from the vector
in (5.10) after this transformation agrees with the random vector in
$(5.10')$ if we write $\e_j=1$ for $j\in M$ and $\e_j=-1$ for
$j\notin M$ in it. These random vectors really agree (not only in
distribution) since for all functions $f\in \Cal F$ both vectors
have a component which is the sum of terms of the form
$f(\xi_{j_1,1}^{(\delta_{j_1})},\dots,\xi_{j_k,k}^{(\delta_{j_k})})$,
$\delta_{j_s}=\pm1$, $1\le s\le k$, multiplied with an appropriate
power of $(-1)$, and this power equals the number of $-1$ components
in the sequence $\delta_{j_1},\dots,\delta_{j_k}$ plus the cardinality
of the set $\{j_1,\dots,j_k\}\cap M$. Part b) of the lemma can be
proved in the same way, hence it is omitted.
 
\medskip
Lemma 4A will be proved with the help of part a) of Lemma~5 and the
following Lemma~6A.
\medskip\noindent
{\bf Lemma 6A.} {\it Let us consider a class of functions $\Cal F$
satisfying the conditions of Proposition 4, and the random variables
$\bar I_{n,k}^V(f)$, $f\in\Cal F$, $V\subset\{1,\dots,k\}$, defined
in formula (5.1). Let $\Cal B=\Cal
B(\xi_{1,s}^{(1)},\dots,\xi_{n,s}^{(1)};\;1\le s\le k)$ denote the
$\sigma$-algebra generated by the random variables
$\xi_{1,s}^{(1)},\dots,\xi_{n,s}^{(1)}$ , $1\le s\le k$, taking part in
the definition of the random variables $\bar I_{n,k}^V(f)$. For all
$V\subset\{1,\dots,k\}$, $V\neq\{1,\dots,k\}$, there exists a number
$A_0=A_0(k)>0$ such that the inequality
$$
P\(\sup_{f\in\Cal F}\left. E\(\bar I_{n,k}^V(f)^2\right|\Cal B\)
> 2^{-(3k+3)}A^2n^{2k}\sigma^{2k+2}\)<
n^{k-1}e^{-A^{1/(2k-1)} n\sigma^2/k}
\tag5.12
$$
holds for all $A\ge A_0$.}
\medskip\noindent
{\it Proof of Lemma 6A.}\/ Let us first consider the case
 $V=\emptyset$. Then $\left.E\(\bar I_{n,k}^\emptyset(f)^2\right|
\Cal B\) =E\(\bar I_{n,k}^\emptyset(f)^2\)\le\frac{n!}{k!}\sigma^2
\le n^{2k}\sigma^{2k+2}$ for all $f\in\Cal F$. In the above
calculation we exploited that the functions $f\in\Cal F$ are
canonical, and this implies certain orthogonalities, and beside
this the inequality $n\sigma^2\ge1$ holds. The above relation
implies inequality (5.12) for $V=\emptyset$ for all $\oo\in\Omega$
if the number $A_0$ is chosen sufficiently large.
 
To avoid some complications in the notation let us restrict our
attention to the sets $V=\{1,\dots,u\}$, $1\le u<k$, and prove
relation (5.12) for such sets. For this goal let us introduce the
random variables
$$
\bar I_{n,k}^V(f,j_{u+1},\dots,j_k)=\frac1{k!}\summ\Sb 1\le j_s\le
n,\; s=1,\dots, u\\ j_s\neq j_{s'} \text{ if } s\neq s',\;1\le
s,s'\le k\endSb
f\(\xi_{j_1,1}^{(1)},\dots,\xi_{j_u,u}^{(1)},\xi_{j_{u+1},u+1}^{(-1)},
\dots, \xi_{j_k,k}^{(-1)}\),
$$
for all $f\in\Cal F$, i.e. we fix some indices $j_{u+1},\dots,j_k$,
$1\le j_s\le n$, $u+1\le s\le k$, $j_s\neq j_{s'}$ if $s\neq s'$, and
sum up only those terms in the sum defining $\bar I_{n,k}^V(f)$
which contain $\xi_{j_{u+1},u+1}^{(-1)},\dots, \xi_{j_k,k}^{(-1)}$ in
their last $k-u$ coordinates. Then we can write
$$
\aligned
\left.E\(\bar I_{n,k}^V(f)^2\right|\Cal B\)&=
\left.E\(\(\summ\Sb 1\le j_s\le n\; s=u+1,\dots,k\\ j_s\neq j_{s'}\text
{if } s\neq s'\endSb
\bar I_{n,k}^V(f,j_{u+1},\dots,j_{u_k})\)^2\right|\Cal B\) \\
&=\summ\Sb 1\le j_s\le n\; s=u+1,\dots,k\\ j_s\neq j_{s'}\text
{if } s\neq s'\endSb
\left.E\(\bar I_{n,k}^V(f,j_{u+1},\dots,j_{u_k})^2\right|\Cal B\).
\endaligned \tag5.13
$$
The last relation follows from the identity
$$
\left.E\(\bar I_{n,k}^V(f,j_{u+1},\dots,j_{u_k})
\bar I_{n,k}^V(f,j'_{u+1},\dots,j'_{u_k})\right|\Cal B\)=0
$$
if $(j_{u+1},\dots,j_k)\neq(j'_{u+1},\dots,j'_k)$, which relation
holds, since $f$ is a canonical function.
 
It follows from relation (5.13) that
$$
\aligned
&\left\{\oo\:\sup_{f\in\Cal F}\left.E\(\bar I_{n,k}^V(f)^2\right|
\Cal B\)(\oo) > 2^{-(3k+3)}A^2n^{2k}\sigma^{2k+2}  \right\}\\
&\qquad \subset \bigcup \Sb 1\le j_s\le n\; s=u+1,\dots,k\\
j_s\neq j_{s'} \text {if } s\neq s'\endSb
\left\{\oo\: \sup_{f\in\Cal F}\left. E\(\bar
I_{n,k}^V(f,j_{u+1},\dots,j_{u_k})^2\right|\Cal
B\)(\oo)>\frac{A^2n^{2k}\sigma^{2k+2}}{2^{(3k+3)}n^{k-u}} \right\}.
\endaligned
\tag5.14
$$
The probability of the events in the union at the right-hand side
of (5.14) can be estimated with the help of the corollary of
Proposition~5 with parameter $u<k$ instead of $k$. (We may assume
that Proposition~5 holds for $u<k$.) This corollary yields that
$$
P\(\sup_{f\in\Cal F}\left. E\(\bar
I_{n,k}^V(f,j_{u+1},\dots,j_{u_k} )^2\right|\Cal
B\)>\frac {A^2\sigma^{2k+2}n^{k+u}} {2^{(2k+3)}}\)\le
e^{-A^{-1/(2u+1)}(n-u)\sigma^2}. \tag5.15
$$
Indeed, the expression $\left. E\(\bar
I_{n,k}^V(f,j_{u+1},\dots,j_{u_k} )^2\right|\Cal B\)$
can be calculated in the following way: Take the independent
$U$-statistic
$$
\bar I_{n,k}^V(f,x_{u+1},\dots,x_k)=\frac1{k!}\summ\Sb
j_s\in\{1,\dots,n\}\setminus\{j_{u+1},\dots,j_k\},\\
s=1,\dots, u,\; j_s\neq j_{s'} \text{ if } s\neq s'\endSb
f\(\xi_{j_1,1}^{(1)},\dots,\xi_{j_u,u}^{(1)},x_{u+1},\dots,x_k\),
\tag5.16
$$
of order $u$ with sample size $n-k+u$, and integrate the square of
this function with respect to the variables $x_{u+1},\dots,x_k$ by
the measure $\mu^{k-u}$. Hence the expression  at the left-hand side
of (5.15) can be bounded by means of Proposition 5 if we apply it
for our class of functions $\Cal F$  considering them as functions
on $(X^u\times Y, \Cal X^u\times \Cal Y, \mu^u\times\rho)$ with
$(Y,\Cal Y,\rho)=(X^{k-u},\Cal X^{k-u},\mu^{k-u})$. (A small
inaccuracy was committed in the above statement because to
define the expression in (5.16) as a $U$-statistic we should have
divided by $u!$ instead of $k!$. But this causes no real problem.)
 
We get inequality (5.15) from Proposition~5 by replacing the level
$\frac {A^2\sigma^{2k+2}n^{k+u}} {2^{(3k+3)}}$ in the
probability at the left-hand side by $A^2(n-u)^{2u} \sigma^{2u+2}
<\frac {A^2\sigma^{2k+2}n^{k+u}}{2^{(2k+2)}}$. The last
inequality really holds if the constant $K$ is chosen sufficiently
large in the condition $n\sigma^2>K\log n$ of Proposition~4.
 
Relations (5.14) and (5.15) imply that
$$
P\(\sup_{f\in\Cal F}\left. E\(\bar I_{n,k}^V(f)^2\right|
\Cal B\)(\oo) > 2^{-(3k+3)}A^2n^{2k}\sigma^{2k+2} \)\le
n^{k-u}e^{-A^{-1/(2u+1)}(n-u)\sigma^2},
$$
and $u\le k-1$. Hence also inequality (5.12) holds.
\medskip
Now we prove Lemma~4A.
\medskip\noindent
{\it Proof of Lemma 4A.} We show with the help of Lemmas 2 and
Lemma~6A that
$$
\aligned
P\(\sup_{f\in\Cal F} n^{k/2}\left|\bar
I_{n,k}(f)\right|>An^{k/2}\sigma^{k+1}\)&<
2P\(\sup_{f\in\Cal F} |S(f)|>\frac A2n^k\sigma^{k+1}\)\\
&\qquad +2^kn^{k-1}e^{-A^{1/(2k-1)} n\sigma^2/k}
\endaligned \tag5.17
$$
with the function $S(f)$ defined in (5.10). To prove relation
(5.17) introduce the random variables
$Z(f)=(-1)^{k}\bar I_{n,k}^{\{1,\dots,k\}}(f)$
and $\bar Z(f)=\summ_{V\subset \{1,\dots,k\},\,
V\neq\{1,\dots,k\}}(-1)^{|V|+1}\bar I_{n,k}^V(f)$ for all
$f\in\Cal F$, the $\sigma$-algebra $\Cal B$ considered in Lemma~6A
and the set
$$
B=\bigcap\Sb V\subset\{1,\dots,k\}\\V\neq\{1,\dots,k\}\endSb
\left\{\oo\: \sup_{f\in\Cal F}\left.E\(\bar I_{n,k}^V(f)^2\right|
\Cal B\)(\oo) \le 2^{-(3k+3)}A^2n^{2k}\sigma^{2k+2}\right\}.
$$
 
Observe that $S(f)=Z(f)-\bar Z(f)$, $f\in\Cal F$, $B\in\Cal B$,
and by Lemma~6A  the inequality $1-P(B)\le2^kn^{k-1}
e^{-A^{1/(2k-1)} n\sigma^2/k}$ holds. Hence to prove relation
(5.17) as a consequence of Lemma~2 it is enough to show that
$$
\left.P\(|\bar Z(f)|>\frac A2n^k\sigma^{k+1}\right|\Cal
B\)(\oo)\le\frac12 \quad \text{ for all }f\in\Cal F \quad \text {if }
\oo\in\Cal B. \tag5.18
$$
But $P\(\bar I_{n,k}(f)|>2^{-(k+1)} An^k\sigma^{k+1}|\Cal
F\)(\oo)\le 2^{-(k+1)}$ for all $f\in \Cal F$ if $\oo\in B$ by the
`conditional Chebishev inequality', hence relation (5.18) holds.
 
Lemma 4A follows from relation (5.17), part~a of Lemma~5 and the
observation that the random vectors $\{\bar I_{n,k}^{(V,\e)}(f)\}$,
$f\in\Cal F$, defined in $(5.9')$ have the same distribution for all
$V\subset\{1,\dots,k\}$ as the random vector $\bar I_{n,k}^{\e}(f)$,
$f\in\Cal F$, considered in the formulation of Lemma~4A. Hence
$$
P\(\sup_{f\in\Cal F} |S(f)|>\frac A2n^k\sigma^{k+1}\)\le
2^kP\(\sup_{f\in\Cal F} \left|\bar I_{n,k}^{\e}(f)\right|
>2^{-(k+1)}A n^k\sigma^{k+1}\).
$$
\medskip
 
In the proof of Lemma~4B we apply following Lemma~6B which is a
version of Lemma~6A.
\medskip\noindent
{\bf Lemma 6B.} {\it Let us consider a class of functions $\Cal F$
satisfying the conditions of Proposition~5 and the random variables
$\bar I_{n,k}^V(f,y)$, $f\in\Cal F$, $V\subset\{1,\dots,k\}$,
defined in formula (5.3). Let $\Cal B=\Cal B(\xi_{1,s}^{(1)},
\dots, \xi_{n,s}^{(1)};\;1\le s\le k)$ denote the $\sigma$-algebra
generated by the random variables
$\xi_{1,s}^{(1)},\dots,\xi_{n,s}^{(1)}$, $1\le s\le k$, taking part
in the definition of the random variables $\bar
I_{n,k}^V(f,y)$ and $H_{n,k}^V(f)$.
\medskip
\item{a)} For all $V\subset\{1,\dots,k\}$, $V\neq\{1,\dots,k\}$,
there exists a number $A_0=A_0(k)>0$ such that the inequality
$$
P\(\sup_{f\in\Cal F} E(H^{V}_{n,k}(f)|\Cal B)
>2^{-(4k+4)}A^{(2k-1)/k} n^{2k}\sigma^{2k+2}\)<
n^{k-1}e^{-A^{1/2k} n\sigma^2/k}.
\tag5.19
$$
holds for all $A\ge A_0$.
\medskip
\item{b)} Given two subsets $V_1,V_2\subset\{1,\dots,k\}$ of the
set $\{1,\dots,k\}$ define the random integrals
$$
H_{n,k}^{(V_1,V_2)}(f)=\int |\bar I_{n,k}^{V_1}(f,y)
\bar I_{n,k}^{V_2}(f,y)| \rho(\,dy),
\quad f\in\Cal F,
$$
with the help of the functions $\bar I_{n,k}^V(f,y)$ defined in
(5.3). If at least one of the sets $V_1$ and $V_2$ is not the
set $\{1,\dots,k\}$, then there exists some number $A_0=A_0(k)>0$
such that if the integrals $H_{n,k}(f)$, $f\in\Cal F$, determined by
this class of functions $\Cal F$ have a good tail behaviour at level
$T^{(2k+1)/2k}$ for some $T\ge A_0$, then the inequality
$$
P\(\sup_{f\in\Cal F} E(H^{(V_1,V_2)}_{n,k}(f)|\Cal B)
>2^{-(2k+2)}A^2n^{2k}\sigma^{2k+2}\)<2n^{k-1}e^{-A^{1/2k}n\sigma^2/k}.
\tag5.20
$$
holds for all $A\ge T$.}
\medskip\noindent
{\it Proof of Lemma 6B.}\/ Part a) of Lemma 6B can be proved in the
same way as Lemma 6A, only the formulas applied in the proof become a
little bit more complicated. Hence I omit the proof. (The difference
between the power of the parameter $A$ at the right-hand side of
formulas (5.19) and (5.12) appear, since the left-hand side of
(5.19) contains the term $A^{(2k-1)/2k}$ and not $A^2$.)
Part b) will be proved with the help of Part a) and the inequality
$$
\sup_{f\in\Cal F} E(H^{(V_1,V_2)}_{n,k}(f)|\Cal B) \le
\(\sup_{f\in\Cal F} E(H^{V_1}_{n,k}(f)|\Cal B)\)^{1/2}
\(\sup_{f\in\Cal F} E(H^{V_2}_{n,k}(f)|\Cal B)\)^{1/2}
$$
which follows from the Schwarz inequality applied for integrals with
respect to conditional distributions. Let us assume that
$V_1\neq\{1,\dots,k\}$. The last inequality implies that
$$
\aligned
&P\(\sup_{f\in\Cal F} E(H^{(V_1,V_2)}_{n,k}(f)|\Cal B)
>2^{-(2k+2)}A^2n^{2k}\sigma^{2k+2}\)\\
&\qquad \le P\(\sup_{f\in\Cal F} E(H^{V_1}_{n,k}(f)|\Cal B)
>2^{-(4k+4)}A^{(2k-1)/k} n^{2k}\sigma^{2k+2}\) \\
&\qquad\qquad+P\(\sup_{f\in\Cal F} E(H^{V_2}_{n,k}(f)|\Cal B)
>A^{(2k+1)/k} n^{2k}\sigma^{2k+2}\)
\endaligned
$$
Hence the estimate (5.19) for $V=V_1$ together with the inequality
$$
P\(\sup_{f\in\Cal F} E(H^{V_2}_{n,k}(f)|\Cal B)
>A^{(2k+1)/k} n^{2k}\sigma^{2k+2}\)\le n^{k-1} e^{-A^{1/2k}n\sigma^2/k}
$$
which follows from Part a) if $V_2\neq\{1,\dots,n\}$ (in this case
the level $A^{(2k+1)/k} n^{2k}\sigma^{2k+2}$ can be replaced by
$2^{-(4k+4)}A^{(2k-1)/k} n^{2k}\sigma^{2k+2}$ in the probability we
consider) and from the conditions of Part b) if
$V_2=\{1,\dots,k\}$ imply relation (5.20).
\medskip
Now I prove Lemma~4B.
\medskip\noindent
{\it Proof of Lemma 4B.}\/ By Part b) of Lemma~5 it is enough to
prove that relation (5.6) holds if the random variables $\bar W(f)$
are replaced in it by the random variables $W(f)$ defined in
formula~(5.11). We shall prove this by applying Lemma~2 with the
choice of $Z(f)=H_{n,k}^{(\bar V,\bar V)}(f)$, $\bar V=\{1,\dots,k\}$,
$\bar Z(f)=W(f)-Z(f)$, $f\in\Cal F$, $\Cal B=\Cal
B(\xi_{1,s}^{(1)},\dots, \xi_{n,s}^{(1)};\;1\le s\le k)$, and the set
$$
B=\bigcap\Sb (V_1,V_2)\: V_j\subset \{1,\dots,k\},\;j=1,2\\
V_1\neq\{1,\dots,k\} \text { or } V_2\neq\{1,\dots,k\} \endSb
\left\{\oo\: \sup_{f\in\Cal F} E(H^{(V_1,V_2)}_{n,k}(f)|\Cal B)(\oo)
\le 2^{-(2k+2)} A^{2} n^{2k}\sigma^{2k+2}\right\}.
$$
 
By Lemma 6B $1-P(B))\le2^{2k+1}n^{k-1} e^{-A^{1/2k}n\sigma^2/k}$,
and to prove Lemma 4B with the help of Lemma~2 it is enough to show
that
$$
P\(\left.|\bar Z(f)|>\frac{A^2}2 n^{2k}\sigma^{2(k+1)}\right|
\Cal B\)(\oo)\le\frac12 \quad \text{for all }f\in \Cal F \text { if }
\oo\in B.
$$
To prove this relation observe that
$$
E(|\bar Z(f)| |\Cal B)\le \summ \Sb (V_1,V_2)\: V_j\subset
\{1,\dots,k\},\;j=1,2\\
V_1\neq\{1,\dots,k\} \text { or } V_2\neq\{1,\dots,k\} \endSb
E(H^{(V_1,V_2)}_{n,k}(f)|\Cal B), \le\frac{A^2}4n^{2k}\sigma^{2k+2}
\quad \text{if } \oo\in B
$$
for all $f\in \Cal F$. Hence the `conditional Markov inequality'
implies that
$$
P\(\left. |\bar Z(f)|> \frac{A^2}2n^{2k}\sigma^{2k+2}\right|\Cal B\)
\le\frac12 \quad\text{if }\oo\in B\quad \text{and } f\in\Cal F.
$$
Lemma~4B is proved.
 
\beginsection
The proof of Propositions 4 and 5
 
The proof of Propositions 4 and 5 for general $k\ge1$ with the help
of the symmetrization lemmas~4A and~4B is similar to the proof of
Proposition~4 in the case $k=1$ presented in Section~4. The proof
applies an induction procedure with respect to the parameter $k$.
In the proof of Proposition~4 for parameter~$k$ we may assume that
Propositions~4 and~5 hold for $k'<k$. In the proof of Proposition~5
we may also assume that Proposition~4 holds for the parameter~$k$.
 
In the proof of Proposition 4 let us introduce (with the notation of
this proposition) the functions
$$
S^2_{n,k}(f)(x_{j,s},\,1\le j\le n,\,1\le s\le k)=\frac1{k!}\summ\Sb
1\le j_s\le n,\; s=1,\dots, k\\ j_s\neq j_{s'} \text{ if } s\neq
s'\endSb f^2\(x_{j_1,1},\dots,x_{j_k,k}\),\quad f\in\Cal F, \tag6.1
$$
where $x_{j,s}\in X$, $1\le j\le n$, $1\le s\le k$. Fix some number
$A>T$ and define the set $H$
$$
\aligned
H=H(A)&=\biggl\{(x_{j,s},\,1\le j\le n,\,1\le s\le k), \\
&\qquad \sup_{f\in\Cal F} S^2_{n,k}(f)(x_{j,s},\,1\le j\le n,\,
1\le s\le k)>2^kA^{4/3}n^k\sigma^2\biggr\}.
\endaligned \tag6.2
$$
We want to show that
$$
P(\{\oo\:(\xi_{j,s}(\oo),\,1\le j\le n,\,1\le s\le k)\in H\})\le 2^k
e^{-A^{2/3k}n\sigma^2} \quad\text{if }A\ge T. \tag6.3
$$
 
Relation (6.3) will be proved by means of the Hoeffding decomposition
of the $U$-statistics with kernel functions $f^2(x_1,\dots,x_k)$,
$f\in\Cal F$, and by the estimation of the sum this decomposition
yields. More explicitly, write
$$
f^2(x_1,\dots,x_k)=\summ_{V\subset\{1,\dots,k\}} f_V(x_j,j\in V)
\tag6.4
$$
with $f_V(x_j,j\in V)=\prodd_{j\notin V}P_j\prodd_{j\in V}Q_j
f^2(x_1,\dots,x_k)$, where $P_j$ and $Q_j$ are the operators $P_\mu$
and $Q_\mu$ defined in formulas (2.6) and $(2.6'')$ if $(Y_1\times
Z\times Y_2,\Cal Y_1\times \Cal Z\times Y_2)$ is the $k$-fold product
$(X^k,\Cal X^k)$ of the measure space $(X,\Cal X)$ in these
definitions, $Z$ is the $j$-th component in these products, and $Y_1$
is the product of the components before and $Y_2$ is the product of
the components after this component. (Relation (6.4) follows from the
identity $f^2=\prodd_{j=1}^k(P_j+Q_j)f^2$ if the multiplications are
carried out in this formula. In the calculation we exploit that the
operators $P_j$ and $P_{j'}$ are commutative if $j\neq j'$, and
the same relation holds for the pairs $P_j$ and $Q_{j'}$ or $Q_j$ and
$P_{j'}$ or $Q_j$ and $Q_{j'}$.)
 
The identity $S^2_{n,k}(f)(\xi_{j,r}\,1\le j\le n,1\le r\le
k)=k!I_{n,k}(f^2)$ holds for all $f\in\Cal F$, and by writing the
Hoeffding decomposition (6.4) for each term
$f^2(\xi_{j_1,1}\dots,\xi_{j_k,k})$ in the expression
$I_{n,k}(f^2)$ we get that
$$
\aligned
&P\(\sup_{f\in\Cal F}S^2_{n,k}(f)(\xi_{j,s},\,1\le j\le n,\,1\le s\le k)
>2^kA^{4/3}n^k\sigma^2\)\\
&\qquad \le\summ_{V\subset\{1,\dots,k\}} P\(\sup_{f\in\Cal F}
n^{k-|V|}|\bar I_{n,|V|}(f_V)|>A^{4/3}n^k\sigma^2\)
\endaligned \tag6.5
$$
with the functions $f_V$ in (6.4). We want to give a good
estimate for all terms in the sum at the right-hand side in (6.5).
For this goal we show that the classes of functions $\{f_V\: f\in
\Cal F\}$ satisfy the conditions of Proposition 4 for all
$V\subset\{1,\dots,k\}$.
 
The functions $f_V$ are canonical for all $V\subset\{1,\dots,k\}$.
(This follows from the commutativity relations between the
operators $P_j$ and $Q_j$ mentioned before, the identity $P_jQ_j=0$
and the fact that the canonical property of the function can be
expressed in the form $P_jf_V=0$ for all $j\in V$.) We have
$|f^2(x_1,\dots,x_k)|\le 2^{-2(k+1)}$. The norm of $Q_j$ as a map
from the $L_\infty$ space to $L_\infty$ space is less than 2, the
norm of $P_j$ is less than 1, hence $\left|\supp_{x_j\in X,j\in
V}f_V(x_j,j\in V)\right|\le 2^{-(k+2)}\le2^{-(k+1)}$ for all
$V\subset\{1,\dots,k\}$. We have $\int
f^4(x_1,\dots,x_k)\mu(\,dx_1)\dots\mu(\,dx_k)\le 2^{-(k+1)}\sigma^2$,
hence $\int f^2_V(x_j,j\in V)\prodd_{j\in V}\mu(\,dx_j)\le
2^{-(k+1)}\sigma^2\le\sigma^2$ for all $V\subset\{1,\dots,k\}$ by
Lemma~1. Finally, to check that the class of functions $\Cal
F_V=\{f_V\:f\in\Cal F\}$ is $L_2$-dense with exponent $L$ and
parameter $D$ observe that for all probability measures $\rho$ on
$(X^k,\Cal X^k)$ and pairs of functions $f,g\in \Cal F$
$\int(f^2-g^2)^2\,d\rho\le 2^{-2k}\int(f-g)^2\,d\rho$. This implies
that if $\{f_1,\dots,f_m\}$, $m\le D\e^{-L}$, is an $\e$-dense
subset of $\Cal F$ in the space $L_2(X^k,\Cal X^k,\rho)$, then the
set of functions $\{2^kf_1^2,\dots,2^kf_m^2\}$ is an $\e$-dense
subset of the class of functions $\Cal F'=\{2^kf^2\: f\in \Cal F\}$.
Then by Lemma~1 for all $V\subset\{1,\dots,k\}$ the set of functions
$\{(f_1)_V,\dots,(f_m)_V)$ is an $\e$-dense subset of the class of
functions $\Cal F_V$ in the space $L_2(X^k,\Cal X^k,\rho)$. This
means that $\Cal F_V$ is also $L_2$-dense with exponent $L$ and
parameter~$D$.
 
For $V=\emptyset$ the relation $f_V=\int f^2(x_1,\dots,x_k)
\mu(\,dx_1)\dots\mu(\,dx_k)\le\sigma^2$ holds, and
$I_{|V|}(f_{|V|})|=f_V\le\sigma^2$. Therefore the term corresponding
to $V=\emptyset$ in the sum at the right-hand side
of (6.5) equals zero if $A_0\ge1$ in the conditions of Proposition~4.
The terms corresponding to sets $V$, $1\le|V|\le k$ in these sums
satisfy the inequality
$$
\align
&P\(\sup_{f\in\Cal F}|\bar I_{n,|V|}(f_V)|>A^{4/3}n^{|V|}\sigma^2\)\\
&\qquad \le P\(\sup_{f\in\Cal F}
|\bar I_{n,|V|}(f_V)|>A^{4/3}n^{|V|}\sigma^{|V|+1}\)
\le e^{-A^{2/3k}n\sigma^2} \quad\text{if } 1\le|V|\le k.
\endalign
$$
This inequality follows from the inductive hypothesis if $|V|<k$, and
in the case $V=\{1,\dots,k\}$  from the inequality $A\ge
T$ and the assumption that $U$-statistics determined by a class of
functions satisfying the conditions of Proposition~4 have a good
tail behaviour at level $T^{4/3}$. The last relation together with
formula (6.5) imply relation~(6.3).
 
By conditioning the probability $P\(\left|\bar I_{n,k}^{\e}(f)
\right|>2^{-(k+2)}A n^{k/2}\sigma^{k+1}\)$ with respect to the
random variables $\xi_{j,s}$, $1\le j\le n$, $1\le s\le k$ we get
with the help of Proposition~A that
$$
\align
&P\(\left.\left|\bar I_{n,k}^{\e}(f)\right|
>2^{-(k+2)}A n^k\sigma^{k+1}\right|\xi_{j,s}(\oo)=x_{j,s},
1\le j\le n,1\le s\le k\) \\
&\qquad \le C\exp\left\{-B\(\frac{A^2n^{2k}\sigma^{2(k+1)}}{2^{2k+4}
S^2_{n,k}(x_{j,s},1\le j\le n,1\le s\le k)}\)^{1/k}\right\} \tag6.6 \\
&\qquad \le Ce^{-2^{-3-4/k}BA^{2/3k}n\sigma^2} \quad \text{for all }
f\in\Cal F\quad \text{if }\{x_{j,s},\, 1\le j\le n,\,1\le s\le
k\}\notin H.
\endalign
$$
Given some points $x_{j,s}$, $1\le j\le n$, $1\le s\le k$, define the
probability measures $\rho_s$, $1\le s\le k$, uniformly distributed
on the set $x_{j,s}$, $1\le j\le s$, i.e. $\rho_s(x_{j,s})=\frac1n$,
$1\le j\le n$, and their product $\rho=\rho_1\times\cdots\times\rho_k$.
If $f$ is a function on $(X^k,\Cal X^k)$ such that $\int f^2
\,d\rho\le\delta^2$ with some $\delta>0$, then $|f(x_{j,s})|
\le \delta n^{k/2}$ for all $1\le s\le k$, $1\le j\le n$, and
$P\(\left.\left|\bar I_{n,k}^{\e}(f)\right|>\delta n^{3k/2}
\right|\xi_{j,s}=x_{j,s}, 1\le j\le n,\, 1\le s\le k\)=0$. Choose
the numbers $\bar\delta=An^{-k/2}2^{-(k+2)}\sigma^{k+1}$ and
$\delta=2^{-(k+2)}n^{-k-1/2}\le\bar\delta$. (The inequality
$\delta\le\bar\delta$ holds, since $A\ge A_0\ge1$, and $\sigma\ge
n^{-1/2}$.)  Choose a $\delta$-dense set $\{f_1,\dots,f_m\}$ in the
$L_2(X^k,\Cal X^k,\rho)$ space with $m\le D\delta^{-L}\le 2^{(k+2)L}
n^{\beta+(k+1/2)L}$ elements. Then formula (6.6) implies that
$$
\align
&P\(\sup_{f\in\Cal F}\left.\left|\bar I_{n,k}^{\e}(f)\right|
>2^{-(k+1)}A n^k\sigma^{k+1}\right|\xi_{j,s}(\oo)=x_{j,s},
1\le j\le n,1\le s\le k\) \\
&\quad \le \sum_{j=1}^m P\(\left.\left|\bar I_{n,k}^{\e}(f_j)\right|
>2^{-(k+2)}A n^k\sigma^{k+1}\right|\xi_{j,s}(\oo)=x_{j,s},
1\le j\le n,1\le s\le k\)  \tag6.7    \\
&\qquad \le C 2^{(k+2)L}n^{\beta+(k+1/2)L}
e^{-2^{-3-4/k}BA^{2/3k}n\sigma^2} \quad \text{if }\{x_{j,s},\, 1\le
j\le n,\,1\le s\le k\}\notin H.
\endalign
$$
 
Relations (6.3) and (6.7) imply that
$$
\aligned
&P\(\sup_{f\in\Cal F}.\left|\bar I_{n,k}^{\e}(f)\right|
>2^{-(k+1)}A n^k\sigma^{k+1}\) \\
&\qquad \le C2^{(k+2)L}n^{\beta+(k+1/2)L}
e^{-2^{-3-4/k}BA^{2/3k}n\sigma^2}+
2^k e^{-A^{2/3k}n\sigma^2} \quad\text{if }A\ge T.
\endaligned \tag6.8
$$
Proposition 4 follows from the estimates (5.2) and (6.8) if the
constants $A_0$ and $K$ in the condition $n\sigma^2\ge K((L+\beta)
\log n+1)$ are chosen sufficiently large. In this case the upper
bound these estimate yields for the probability at the left-hand side
of (3.8) is smaller than $e^{-A^{2/k}n\sigma^2}$.
\medskip
 
Let us turn to the proof of Proposition~5. By formula (5.8) it is
enough to show that
$$
\aligned
&P\(\sup_{f\in\Cal F} \left |H_{n,k}(f|G,V_1,V_2)\right|
>\frac{A^2}{2^{4k+1}k!} n^{2k}\sigma^{2(k+1)}\) \le
e^{-A^{1/2k}n\sigma^2}\\
&\qquad\text{ for all } G\in \Cal G\quad \text{and }
\;V_1,V_2\in\{1,\dots,k\} \quad\text{if } A\ge A_0.
\endaligned \tag6.9
$$
with the random variables $H_{n,k}(f|G,V_1,V_2)$ defined in formula
(5.7). Let us first prove (6.9) in the case when $|e(G)|=k$, i.e.\
all vertices of the diagram $G$ are an end-point of some edge, and
the expression $H_{n,k}(f|G,V_1,V_2)$ contains no `symmetryzing term'
$\e_j$.  By the Schwarz inequality
$$
\aligned
|H_{n,k}(f|G,V_1,V_2)|&\le
\(\sum\Sb j_1,\dots,j_k, 1\le j_s\le n,\\ j_s\neq j_{s'} \text{ if
} s\neq s'\endSb \int
f^2(\xi_{j_1,1}^{(\delta_1)},\dots,\xi_{j_k,k}^{(\delta_k)},y)
\rho(\,dy)\)^{1/2}\\
&\qquad \(\sum\Sb j_1,\dots,j_k, 1\le j_s\le n,\\ j_s\neq j_{s'}
\text{ if }s\neq s'\endSb \int f^2(\xi_{j_1,1}^{(\bar\delta_1)},
\dots,\xi_{j_k,k}^{(\bar\delta_k)},y) \rho(\,dy)\)^{1/2},
\endaligned \tag6.10
$$
for such diagrams $G$, where $\delta_s=1$ if $s\in V_1$,
$\delta_s=-1$ if $s\notin V_1$, and $\bar\delta_s=1$ if $s\in V_2$,
$\bar\delta_s=-1$ if $s\notin V_2$. Hence
$$
\align
&\left\{\oo\:\sup_{f\in\Cal F} \left |H_{n,k}(f|G,V_1,V_2)(\oo)\right|
>\frac{A^2}{2^{4k+1}k!} n^{2k}\sigma^{2(k+1)}\right\} \\
&\quad \subset
\left\{\oo\:\sup_{f\in\Cal F} \sum\Sb j_1,\dots,j_k, 1\le j_s\le n,\\
j_s\neq j_{s'} \text{ if } s\neq s'\endSb \int
f^2(\xi_{j_1,1}^{(\delta_1)}(\oo),\dots,\xi_{j_k,k}^{(\delta_k)}
(\oo),y) \rho(\,dy)>\frac {A^2n^{2k}\sigma^{2(k+1)}}
{2^{4k+1}k!} \right\}\\
&\qquad \cup
\left\{\oo\:\sup_{f\in\Cal F} \sum\Sb j_1,\dots,j_k, 1\le j_s\le n,\\
j_s\neq j_{s'} \text{ if } s\neq s'\endSb \int
f^2(\xi_{j_1,1}^{(\bar\delta_1)}(\oo),\dots,\xi_{j_k,k}^{(\bar\delta_k)}
(\oo),y)
\rho(\,dy)>\frac{A^2n^{2k}\sigma^{2(k+1)}}{2^{4k+1}k!}.  \right\}
\endalign
$$
 
The last relation implies that
$$
\align
&P\(\sup_{f\in\Cal F} \left |H_{n,k}(f|G,V_1,V_2)\right|
>\frac{A^2}{2^{4k+1}k!} n^{2k}\sigma^{2(k+1)}\) \\
&\qquad \le 2P\(\sup_{f\in\Cal F} \sum\Sb j_1,\dots,j_k,
1\le j_s\le n,\\ j_s\neq j_{s'} \text{ if } s\neq s'\endSb
h_f(\xi_{j_1,1},\dots,\xi_{j_k,k})
>\frac{A^2n^{2k}\sigma^{2(k+1)}}{2^{4k+1}k!}\) \tag6.11
\endalign
$$
with $h_f(x_1,\dots,x_k)=\int f^2(x_1,\dots,x_k,y)\rho(\,dy)$,
$f\in\Cal F$. (In this upper bound we could get rid of the
terms $\delta_j$ and $\bar\delta_j$, i.e. on the dependence of the
expression $H_{n,k}(f|G,V_1,V_2)$ on the sets $V_1$ and $V_2$, since
the probability of the events in the previous formula do not depend
on these terms.)
 
I claim that
$$
P\(\supp_{f\in\Cal F} |\bar I_{n,k}(h_f)|\ge An^k \sigma^2\)\le
2^k e^{-A^{1/2k}n\sigma^2} \quad \text{for }A\ge A_0  \tag6.12
$$
if the constant $A_0$ and $K$ are chosen sufficiently large in
Proposition~5. Relation (6.12) together with the relation
$\frac{n^{2k}\sigma^{2(k+1)}}{2^{4k+1}k!}\ge n^k\sigma^2$ imply
 that the probability at the right-hand side of (6.11) can be
bounded by $2^{k+1}e^{-A^{1/k}n\sigma^2}$, and the estimate
(6.9) holds in the case $|e(G)|=k$. Relation (6.12)
can be proved similarly to formula (6.3) in the proof of
Proposition~4. It is not difficult to check that $0\le\int
h_f(x_1,\dots,x_k) \mu(\,dx_1)\dots\mu(\,dx_k)\le\sigma^2$,
$\sup|h_f(x_1,\dots,x_k)|\le 2^{-2(k+1)}$, and the class of functions
$\Cal H=\{2^kh_f,\; f\in\Cal F\}$ is an $L_2$-dense class with
exponent $L$ and parameter $D$. This means that by applying the
Hoeffding decomposition of the functions $h_f$, $f\in \Cal F$,
similarly to formula (6.4) we get such sets of functions $(h_f)_V$,
$f\in\Cal F$ for all $V\subset \{1,\dots,k\}$ which satisfy
the conditions of Proposition~4. Hence a natural adaptation of the
estimate given for the expression at the right-hand side of (6.5)
yields the proof of formula (6.12). Let us observe that by our
inductive hypothesis the result of Proposition~4 holds also for $k$,
and this allows us to carry out the estimates we need
also for the class of functions $(h_f)_V$, $f\in\Cal F$, with
$V=\{1,\dots,k\}$ if $A\ge A_0$.
 
In the case $e(G)<k$ formula (6.9) will be proved with the help of
Proposition~A. To carry out this proof first an appropriate expression
will be introduced and bounded for all sets $V_1,V_2\subset
\{1,\dots,k\}$ and diagrams $G\in \Cal G$ such that $|e(G)|<k$. To
define the  expression we shall bound first some notations will be
introduced.
 
Let us consider the set $J_0(G)=J_0(G,k,n)$,
$$
\align
J_0(G)&=\{(j_1,\dots,j_k,j'_1,\dots,j'_k)\:1\le j_s,j'_s\le n,\, 1\le
s,s'\le k,\, j_s\neq j_{s'}\text { if } s\neq s', \\
&\qquad j'_s\neq j'_{s'}\text{ if }s\neq s',\, j_s=j'_{s'}\text{ if }
(s,s')\in G\, j_s\neq j'_{s'}\text{ if } (s,s')\notin G\},
\endalign
$$
the set of those sequences $(j_1,\dots,j_k,j'_1,\dots,j'_k)$
which appear as indices in the summation in formula (5.7). I give
a partition of $J_0(G)$ appropriate for our purposes.
 
For this aim let us first define the sets
$M_1=M_1(G)=\{s(1),\dots,s(k-|e(G)|)\}=\{1,\dots,k\}\setminus
v_1(G)$, $s(1)<\cdots<s(k-|e(G)|)$, and
$M_2=M_2(G)=\{\bar s(1),\dots,\bar s(k-|e(G|)\}=\{1,\dots,k\}\setminus
v_2(G)$, $\bar s(1)<\cdots<\bar s(k-|e(G|)$, the sets of those vertices
of the first and second row of the diagram $G$ in increasing order
from which no edges start. Let us also introduce the set
$V(G)=V(G,n,k)$,
$$
\align
V(G)&=\{(j_{s(1)},\dots,j_{s(k-|e(G)|)},
j'_{\bar s(1)},\dots,j'_{\bar s(k-|e(G)|)})\:1\le j_{s(p)},
j'_{\bar s(p)}\le n,\\
&\quad 1\le p\le k-|e(G)|,\, j_{s(p)}\neq j_{s(p')},\,
j'_{\bar s(p)}\neq j'_{\bar s(p')} \text { if }p\neq p',\,
1\le p,p'\le k-|e(G)|\},
\endalign
$$
which is the set consisting from the restriction of the coordinates of
the vectors
$$
(j_1,\dots,,j_k,j'_1,\dots,j'_k)\in J_0(G)
$$
to $M_1\cup M_2$. Given a vector $v\in V(G)$ let $v(r)$,
$1\le r\le k-|e(G)|$, and $\bar v(r)$, $1\le r\le k-|e(G)|$, denote
its coordinates corresponding to the set $M_1$ and $M_2$ respectively.
Put
$$
\align
E_G(v)&=\{(j_1,\dots,j_k,j'_1,\dots,j'_k)\:1\le j_s\le n, \text{ if }
s\in v_1(G), \, 1\le j'_s\le n\text { if } s\in v_2(G),\\
&\qquad j_s\neq j_{s'},\, j'_s\neq j'_{s'}\text{ if }s\neq s', \;
j_s=j'_{s'}\text{ if } (s,s')\in G\text { and } j_s\neq j'_{s'}\text{ if
} (s,s')\notin G\\
&\qquad j_{s(r)}=v(r),\, j'_{\bar s(r)}=\bar v(r),\, 1\le r\le
k-|e(G)|\},\quad v\in V(G),
\endalign
$$
where $\{s(1),\dots,s(k-|e(G)|)\}=M_1$, $\{\bar s(1),\dots,
\bar s(k-|e(G)|)\}=M_2$ in the last line of this definition.  The set
$E_G(v)$ contains those vectors in $J_0(G)$ whose coordinates in
$M_1\cup M_2$ are prescribed by the vector $v\in V(G)$ and the
remaining coordinates are chosen freely.
 
Now we define the partition
$$
J_0(G)=\bigcup_{v\in V(G)} E_G(v).
$$
of the set $J_0(G)$.
 
The inequality
$$
P\(S(\Cal F|G,V_1,V_2))>A^{8/3}n^{2k}\sigma^4\)\le
2^{k+1}e^{-A^{2/3k}n\sigma^2} \quad\text{if }A\ge A_0\text{ and }
e(G)<k \tag6.13
$$
will be proved for the random variable
$$
\aligned
S(\Cal F|G,V_1,V_2)=\sup_{f\in\Cal F}\sum_{v\in V(G)}
\biggl(&\sum_{(j_1,\dots,j_k,j'_1,\dots,j'_k)\in E_G(v)}
\int f(\xi_{j_1,1}^{(\delta_1)},\dots,\xi_{j_k,k}^{(\delta_k)},y) \\
&\qquad f(\xi_{j'_1,1}^{(\bar\delta_1)},\dots,\xi_{j'_k,k}
^{(\bar\delta_k)},y) \rho(\,dy)\biggr)^2,
\endaligned \tag$6.13'$
$$
where $\delta_s=1$ if $s\in V_1$, $\delta_s=-1$ if $s\notin V_1$,
and $\bar\delta_s=1$ if $s\in V_2$, $\bar\delta_s=-1$ if $s\notin V_2$.

To prove formula (6.13) let us first fix some $v\in V(G)$ and apply
the Schwarz inequality. It yields that
$$
\align
&\(\sum_{(j_1,\dots,j_k,j'_1,\dots,j'_k)\in E_G(v)}
\int f(\xi_{j_1,1}^{(\delta_1)},\dots,\xi_{j_k,k}^{(\delta_k)},y)
f(\xi_{j'_1,1}^{(\bar\delta_1)},\dots,\xi_{j'_k,k}
^{(\bar\delta_k)},y) \rho(\,dy)\)^2\\
&\qquad\le
\(\sum_{(j_1,\dots,j_k,j'_1,\dots,j'_k)\in E_G(v)}
\int f^2(\xi_{j_1,1}^{(\delta_1)},\dots,\xi_{j_k,k}^{(\delta_k)},y)
\rho(\,dy)\) \\
&\qquad\qquad \(\sum_{(j_1,\dots,j_k,j'_1,\dots,j'_k)\in E_G(v)}
f^2(\xi_{j'_1,1}^{(\bar\delta_1)},\dots,\xi_{j'_k,k}
^{(\bar\delta_k)},y) \rho(\,dy)\).
\endalign
$$
for all $v\in V(G)$. Summing up these inequalities for all
$v\in V(G)$ we get that
$$
\align
S(\Cal F|G,V_1,V_2)&\le\sup_{f\in\Cal F}\sum_{v\in V(G)}
\(\sum_{(j_1,\dots,j_k,j'_1,\dots,j'_k)\in E_G(v)}
\int f^2(\xi_{j_1,1}^{(\delta_1)},\dots,\xi_{j_k,k}^{(\delta_k)},y)
\rho(\,dy)\) \\
&\qquad\(\sum_{(j_1,\dots,j_k,j'_1,\dots,j'_k)\in E_G(v)}
f^2(\xi_{j'_1,1}^{(\bar\delta_1)},\dots,\xi_{j'_k,k}
^{(\bar\delta_k)},y) \rho(\,dy)\)   \tag6.14     \\
&\le \supp_{f\in\Cal F} \(\sum_{(j_1,\dots,j_k,j'_1,\dots,j'_k)\in
J_0(G)} \int f^2(\xi_{j_1,1}^{(\delta_1)},\dots,
\xi_{j_k,k}^{(\delta_k)},y) \rho(\,dy)\) \\
&\qquad \qquad \supp_{f\in\Cal F}
\(\sum_{(j_1,\dots,j_k,j'_1,\dots,j'_k)\in J_0(G)}
f^2(\xi_{j'_1,1}^{(\bar\delta_1)},\dots,\xi_{j'_k,k}
^{(\bar\delta_k)},y) \rho(\,dy)\)
\endalign
$$
To check formula (6.14) we have to observe that by multiplying the
inner sum at the left-hand side of this inequality each term
$f^2(\xi_{j_1,1}^{(\delta_1)},\dots,\xi_{j_k,k}^{(\delta_k)},y)
f^2(\xi_{j'_1,1}^{(\bar\delta_1)},\dots,\xi_{j'_k,k}
^{(\bar\delta_k)},y)$ appears only once. (In particular, it is
determined which index $v\in V(G)$ has to be taken in the outer sum
to get this term. The coordinates of this vector $v$ agree with the
coordinates of the vector $j=(j_1,\dots,j_k,j'_1,\dots,j'_k)$ in
$M_1\cup M_2$, with the coordinates of the vector $j$ which
correspond to those vertices from which no edges of the diagram $G$
start.) Beside this, all these products appear if the multiplications
at the right-hand expression are carried out.
 
Relation (6.14) implies that
$$
P(S(\Cal F|G,V_1,V_2))>A^{8/3}n^{2k}\sigma^4) \le
2P\(\supp_{f\in\Cal F} \bar I_{n,k}(h_f)>A^{4/3}n^k\sigma^2\)
$$
with $h_f(x_1,\dots,x_k)=\int f^2(x_1,\dots,x_k,y)\rho(\,dy)$.
(Here we exploited that in the last formula $S(\Cal F|G,V_1,V_2)$
is bounded by the product of two random variables whose distributions
do not depend on the sets $V_1$ and $V_2$.) Thus to prove inequality
(6.13) it is enough to show that
$$
2P\(\supp_{f\in\Cal F} \bar I_{n,k}(h_f)>A^{4/3}n^k\sigma^2\)\le
2^{k+1}e^{-A^{2/3k}} \quad \text{if } A\ge A_0. \tag6.15
$$
Actually formula (6.15) has been already proved, only formula (6.12)
has to be applied, and  the parameter $A$ has to be replaced by
$A^{4/3}$ in it.
 
The proof of Proposition~5 can be completed similarly to
Proposition~4. It follows from Proposition~A that
$$
\aligned
&P\(\left.|H_{n,k}(f|G,V_1,V_2)|
>\frac{A^2}{2^{4k+2}k!} n^{2k}\sigma^{2(k+1)}
\right| \xi^{\pm1}_{j,s},\,1\le j\le n,\,1\le s\le k\)(\oo)\\
&\qquad \le Ce^{-B2^{-(4+2/k)}(k!)^{-1/k} A^{2/3k}n\sigma^2} \quad
\text{if}\quad S(\Cal F|G,V_1,V_2))(\oo)\le A^{8/3}n^{2k}\sigma^4 \\
&\qquad\qquad\text{ for all } f\in\Cal F,\; G\in \Cal G,\;
|e(G)|<k, \quad \text{and } \;V_1,V_2\in\{1,\dots,k\}
\quad\text{if } A\ge A_0.
\endaligned \tag6.16
$$
Indeed, in this case the conditional probability considered in (6.16)
can be bounded by $C\exp\left\{-B\(\frac{A^4n^{4k}\sigma^{4(k+1)}}
{2^{8k+4}(k!)^2S(\Cal F|G,V_1,V_2)}\)^{1/2j}\right\}\le C\exp
\left\{-B\(\frac{A^{4/3}n^{2k}\sigma^{4k}}{2^{8k+4}(k!)^2}\)^{1/2j}
\right\}$, where $2j=2k-2|e(G)|$, the number of vertices of the
diagram $G$ from which no edges start. Since $j\le k$, $n\sigma^2\ge1$,
and also $\frac{A^{4/3}}{2^{8k+4}(k!)^2}\ge1$ if $A_0$ is chosen
sufficiently large the above calculation implies relation~(6.16).
 
Let us show that also the inequality
$$
\aligned
&P\(\left.\sup_{f\in\Cal F} |H_{n,k}(f|G,V_1,V_2)|
>\frac{A^2}{2^{4k+1}k!} n^{2k}\sigma^{2(k+1)}
\right| \xi^{\pm1}_{j,s},\,1\le j\le n,\,1\le s\le k\)(\oo)\\
&\qquad \le Cn^{(3k+1)L/2+\beta}
e^{-BA^{2/3k}n\sigma^2/2^{(4+2/k)}(k!)^{1/k}} \quad
\text{if } S(\Cal F|G,V_1,V_2))(\oo)\le A^{8/3}n^{2k}\sigma^4 \\
&\qquad \qquad\text{ for all } G\in \Cal G,\; |e(G)|<1, \quad
\text{and } \;V_1,V_2\in\{1,\dots,k\} \quad\text{if } A\ge A_0
\endaligned \tag6.17
$$
holds.
 
To deduce formula (6.17) let us fix an elementary event
$\oo\in\Omega$ which satisfies the relation
$S(\Cal F|G,V_1,V_2))(\oo)\le A^{8/3}n^{2k}\sigma^4$,
two sets $V_1,V_2\subset\{1,\dots,k\}$, a
diagram $G$, consider the points $x_{j,s}^{(\pm1)}=
x_{j,s}^{(\pm1)}(\oo)=\xi_{j,s}^{(\pm1)}(\oo)$,
$1\le j\le n$, $1\le s\le k$, and introduce with their help the
following probability measures: For all $1\le s\le k$ define the
probability measures $\nu_s^{(1)}$ which are
uniformly distributed on the points $x_{j,s}^{(\delta_s)}$, $1\le j\le
n$, and $\nu_s^{(2)}$ which are uniformly distributed
on the points $x_{j,s}^{(\bar\delta_s)}$, $1\le j\le n$, where
$\delta_s=1$ if $s\in V_1$, $\delta_s=-1$ if $s\notin V_1$, and
similarly $\bar\delta_s=1$ if $s\in V_2$ and $\bar\delta_s=-1$ if
$s\notin V_2$. Let us consider the product measures
$\alpha_1=\nu_1^{(1)}\times\cdots\times\nu_k^{(1)}\times\rho$,
$\alpha_2=\nu_1^{(2)}\times\cdots\times\nu_k^{(2)}\times\rho$ on
the product space $(X^k\times Y,\Cal X^k\times\Cal Y)$, where $\rho$
is that probability measure on $(Y,\Cal Y)$ which appears in
Proposition~5, together with the measure
$\alpha=\frac{\alpha_1+\alpha_2}2$. Given two functions $f\in \Cal F$
and $g\in\Cal F$ we  give an upper bound for
$|H_{n,k}(f|G,V_1,V_2)(\oo)-H_{n,k}(g|G,V_1,V_2)(\oo)|$ if $\int
(f-g)^2\,d\alpha\le\delta$ with some $\delta>0$. (This bound
does not depend on the `randomizing terms' $\e_j(\oo)$ in the
definition of the random variable $H_{n,k}(\cdot|G,V_1,V_2)$.)
 
In this case $\int(f-g)^2\,d\alpha_j\le2\delta^2$, and
$$
\align
\int&|f(x_{1,j_1}^{(\delta_1)},\dots,x_{k,j_k}^{(\delta_k)},y)-
g(x_{1,j_1}^{(\delta_1)},\dots,x_{k,j_k}^{(\delta_k)},y)|^2
\rho(\,dy) \le2\delta^2n^k, \\
\int& |f(x_{1,j_1}^{(\delta_1)},\dots,x_{k,j_k}^{(\delta_k)},y)-
g(x_{1,j_1}^{(\delta_1)},\dots,x_{k,j_k}^{(\delta_k)},y)|
\rho(\,dy) \le\sqrt2\delta n^{k/2}
\endalign
$$
for all $1\le s\le k$, and $1\le j_s\le n$, and the same result
holds if all $\delta_s$ is replaced by $\bar\delta_s$, $1\le s\le
k$. Since $|f|\le1$ for $f\in\Cal F$, the condition
$\int(f-g)^2\,d\alpha\le \delta^2$ implies that
$$
\align
\int &|f(\xi_{j_1,1}^{(\delta_1)}(\oo),\dots,
\xi_{j_k,k}^{(\delta_k)}(\oo),y)
f(\xi_{j'_1,1}^{(\bar\delta_1)}(\oo),\dots,
\xi_{j'_k,k}^{(\bar\delta_k)}(\oo),y) \rho(\,dy)\\
&\qquad -g(\xi_{j_1,1}^{(\delta_1)}(\oo),\dots,
\xi_{j_k,k}^{(\delta_k)}(\oo),y)
g(\xi_{j'_1,1}^{(\bar\delta_1)}(\oo),\dots,
\xi_{j'_k,k}^{(\bar\delta_k)}(\oo),y) \rho(\,dy)|
\le2\sqrt2\delta n^{k/2}
\endalign
$$
for all vectors $(j_1,\dots,j_k,j'_1,\dots,j'_k)$ which appear as an
index in the summation in (5.7), and
$$
|H_{n,k}(f|G,V_1,V_2)(\oo)-H_{n,k}(g|G,V_1,V_2)(\oo)|
\le2\sqrt2\delta n^{5k/2}
$$
if the originally fixed $\oo\in\Omega$ is considered.
 
Put $\bar\delta=\frac{A^2 n^{-k/2}\sigma^{2(k+1)}}{2^{(4k+7/2)} k!}$,
and $\delta=n^{-(3k+1)/2}\le\bar\delta$ (since $\sigma\ge
\frac1{\sqrt n}$ and we may assume that $A\ge A_0$ is sufficiently
large), choose a $\delta$-dense subset $\{f_1,\dots,f_m\}$ in the
$L_2(X^k\times Y,\Cal X^k\times Y,\alpha)$ space with $m\le
D\delta^{-L}\le n^{(3k+1)L/2+\beta}$ elements. Relation (6.16) for
these functions together with the above estimates yield formula (6.17).
 
It follows from relations (6.13) and (6.17) that
$$
\align
&P\(\sup_{f\in\Cal F}|H_{n,k}(f|G,V_1,V_2)|
>\frac{A^2}{2^{4k+1}k!} n^{2k}\sigma^{2(k+1)}\)\le
2^{k+1}e^{-A^{2/3k}n\sigma^2}\\
&\qquad + Cn^{(3k+1)L/2+\beta}
e^{-BA^{2/3k}n\sigma^2/2^{(4+2/k)}(k!)^{1/k}}
\quad\text{if }A\ge A_0
\endalign
$$
for all $V_1,V_2\subset\{1,\dots,k\}$ also in the case $|e(G)|\le k-1$.
This means that relation (6.9) holds also in this case if the constants
$A_0$ and $K$ are chosen sufficiently large in Proposition~5.
Proposition~5 is proved.
 
\beginsection Appendix. \ The proof of Proposition A
 
The proof will be based on the hypercontractive inequality for
Rademacher functions. Let me first recall this result.
\medskip\noindent
{\bf The hypercontractive inequality for Rademacher functions.}
{\it Let us consider the measure spaces $(X,\Cal X,\mu)$ and $(Y,
\Cal Y,\nu)=(X,\Cal X,\mu)$ defined as $X=\{-1,1\}$, $\Cal X$
contains all subsets of $X$, and $\mu(\{1\})=\mu(\{-1\})=\frac12$.
Given a real number $\gamma>0$ let us define the linear operator
$\bold T_\gamma$ which maps the real (or complex) valued functions on
the space $X$ to the real (or complex) valued functions on the space
$Y$, and satisfies the relations $\bold T_\gamma r_0=r_0$, and
$\bold T_\gamma r_1 =\gamma r_1$, where $r_0(1)=r_0(-1)=1$, and
$r_1(1)=1$, $r_1(-1)=1$. For all $n=1,2,\dots$ let us consider the
$n$-fold product $(X^n, \Cal X^n,\mu^n)$ and $(Y^n,\Cal Y^n,\nu^n)$
together with the $n$-fold product of the operator $\bold T^n_\gamma$
(i.e. $\bold T^n_\gamma$ is the linear transformation for which $\bold
T^n_\gamma (f_1(x_1)\cdots f_n(x_n))=(\bold T_\gamma f_1(x_1)\cdots
\bold T_\gamma f_n(x_n))$ for all functions $f_s$, $1\le s\le n$, on
the space $(X,\Cal X,\mu)$). For all $n=1,2,\dots$ the transformation
$\bold T^n_\gamma$ from the space $L_p(X^n,\Cal X^n,\mu^n)$ to the
space $L_q(Y^n,\Cal Y^n,\nu^n)$ has the norm 1 if $1<p\le q<\infty$,
and $\gamma\le\sqrt{\frac{p-1}{q-1}}$.} \medskip
\medskip
The following corollary of the hypercontractive inequality is useful
for us.
\medskip\noindent
{\bf Corollary of the hypercontractive inequality.}
{\it Let $\e_1,\dots,\e_n$ be independent identically distributed
random variables $P(\e_j=1)=P(\e_j=-1)=\frac12$, $1\le j\le n$,
fix some real numbers $a(j_1,\dots,j_k)$ for all indices
$(j_1,\dots,j_k)$ such that $1\le j_s\le n$, $1\le s\le k$, and
$j_s\neq j_{s'}$ if $s\neq s'$, and define the random variable
$$
Z=\sum\Sb 1\le j_s\le n,\, 1\le s\le k\\ j_s\neq j_{s'}\text{ if }s\neq
s'\endSb a(j_1,\dots,j_k)\e_{j_1}\cdots\e_{j_k}.
$$
The inequality
$$
E|Z|^q\le\(\frac{q-1}{p-1}\)^{kq/2} \(E|Z|^p\)^{q/p}\quad \text{ if }
\quad 1<p\le q<\infty
$$
holds.}
\medskip\noindent
{\it Proof of the Corollary.}\/ Let us define the function
$$
f(x_1,\dots,x_n)=\sum\Sb 1\le j_s\le n,\,
1\le s\le k\\ j_s\neq j_{s'}\text{ if }s\neq s'\endSb
a(j_1,\dots,j_k) r_1(x_{j_1})\cdots r_1(x_{j_k})
$$
on the space $(X^n,\Cal X^n,\mu^n)$. Observe that $\bold
T_\gamma^nf=\gamma^k f$ for this function $f$ and all $\gamma>0$, and
$E|Z|^p=\|f\|_p^p$,
$E|Z|^q=\|f\|_q^q$. Fix some numbers $1<p\le q\le\infty$, and put
$\gamma=\sqrt{\frac{p-1}{q-1}}$. Since the norm of $\bold T^n_\gamma$
as a transformation from the space $L_p(X^n,\Cal X^n,\mu^n)$ to the
space $L_q(Y^n,\Cal Y^n,\nu^n)$ is 1, $\|\bold T_\gamma^n f\|_q
=\gamma^k\|f\|_q\le \|f\|_p$. The above relations imply that
$(E|Z|^q)^{1/q}\le \(\frac{q-1}{p-1}\)^{k/2}E|Z|^p)^{1/p}$ in this
case, and this is what we had to show.
\medskip
Applying the corollary with $p=2$ and some $q>2$ we get that
$$
E|Z|^q\le(q-1)^{kq/2} \(EZ^2\)^{q/2} \le q^{kq/2} \(EZ^2\)^{q/2}=
q^{kq/2} \bar S^q
$$
with
$$
\bar S^2=\sum_{1\le j_1<j_2\cdots<j_k\le n}\(\sum_{\pi\in\Pi_k}
a((j_{\pi(1)},\dots,j_{\pi(k)})\)^2,
$$
where $\Pi_k$ denotes the set of all permutations of the set
$\{1,\dots,k\}$. Observe that
$$
\(\sum_{\pi\in\Pi_k}a(j_{\pi(1)},\dots,j_{\pi(k)})\)^2\le k!
\sum_{\pi\in\Pi_k}
a^2(j_{\pi(1)},\dots,j_{\pi(k)})\quad \text{for all } 1\le j_1<\cdots
j_k\le n,
$$
hence $\bar S^2\le k!S^2$, and $E|Z|^q\le q^{kq/2} (k!)^{q/2}S^q$
with the number $S^2$ defined in (3.4). Thus the Markov inequality
implies that
$$
P(|Z|>x)\le \(q^{k/2}\frac {\sqrt{k!}S}x\)^q   \quad \text{for all
}x>0\quad \text{and } q\ge2.
$$
Choose the number $q$ as the solution of the equation $q\(\frac
{\sqrt{k!}S}x\)^{2/k}=\frac1e$. Then we get
that $P(|Z|>x)\le \exp\left\{- B\(\frac xS\)^{2/k}\right\}$ with
$B=\frac k{2e(k!)^{1/k}}$, provided that $q=\frac1{e{k!}^{1/k}}
\(\frac xS\)^{2/k}\ge2$, i.e. $B\(\frac xS\)^{2/k}\ge k$. By
multiplying the above upper bound with $C=e^k$ we get such an
estimate for $P(|Z|>x)$ which holds for all $x>0$.

\beginsection References:
 
\item{1.)} Alexander, K. (1987) The central limit theorem for empirical
processes over Vapnik--\v{C}ervonenkis classes. {\it Ann. Probab.}
{\bf 15}, 178--203
\item{2.)} Arcones, M. A. and Gin\'e, E. (1994) $U$-processes
indexed by Vapnik--\v{C}ervonenkis classes of functions with
application to asymptotics and bootstrap of $U$-statistics with
estimated parameters. {\it Stoch. Proc. Appl.}  {\bf 52}, 17--38
\item{3.)} Beckner, W. (1975) Inequalities in Fourier Analysis.
Ann. Math. {\bf 102}, 159--182
\item{4.)} de la Pe\~na, V. H. and Gin\'e, E.  (1999)  Decoupling From
dependence to independence. Springer series in statistics. Probability
and its application. Springer Verlag, New York, Berlin, Heidelberg
\item{5.)} de la Pe\~na, V. H. and  Montgomery--Smith, S. (1995)
Decoupling inequalities for the tail-probabilities of multivariate
$U$-statistics. {\it Ann. Probab.}, {\bf 25}, 806--816
\item{6.)} Gross, L. (1975) Logarithmic Soboliev inequalities.
Amer. J. Math.  {\bf 97}, 1061--1083
\item{7.)} Major, P. (1988) On the tail behaviour of the distribution
function of multiple stochastic integrals. {\it Probability Theory
and Related Fields}, {\bf 78},  419--435
\item{8.)} Major, P. (2003) An estimate about multiple stochastic
integrals with respect to a normalized empirical measure.
submitted to {\it Probability Theory and Related Fields}
\item{9.)} Pollard, D. (1984) Convergence of Stochastic Processes.
Springer--Verlag, New York
 
\bigskip\bigskip\bigskip
Supported by the OTKA foundation Nr. 037886
 
\bye